\documentclass[11pt,reqno]{amsart}
\usepackage[utf8]{inputenc}
\usepackage[T1]{fontenc}
\usepackage{lmodern}
\usepackage{amssymb,mathrsfs}
\usepackage{graphicx}
\usepackage[margin=1in,headsep=0.25in]{geometry}
\usepackage[colorlinks,linkcolor=blue,citecolor=red,urlcolor=blue]{hyperref}

\usepackage{tikz}
\usetikzlibrary{arrows.meta,decorations.pathreplacing}

\theoremstyle{plain}
\newtheorem{theorem}{Theorem}[section]
\newtheorem{lemma}[theorem]{Lemma}
\newtheorem{proposition}[theorem]{Proposition}
\newtheorem{corollary}[theorem]{Corollary}
\theoremstyle{remark}

\numberwithin{equation}{section}

\newcommand{\PP}{\mathbb P}
\newcommand{\EE}{\mathbb E}
\newcommand{\ZZ}{\mathbb Z}
\newcommand{\RR}{\mathbb R}

\newcommand{\FF}{\mathcal F}
\newcommand{\GG}{\mathcal G}

\newcommand{\ind}{\mathbf 1}

\author[Yuval Peres]{Yuval Peres}
\address[Yuval Peres]{Beijing Institute of Mathematical Sciences and Applications}
\email{yperes@gmail.com}

\author[Shuo Qin]{Shuo Qin}
\address[Shuo Qin]{Beijing Institute of Mathematical Sciences and Applications, and Yau Mathematical Sciences Center, Tsinghua University}
\email{qinshuo@bimsa.cn}

\makeatletter
\def\@settitle{\begin{center}%
  \baselineskip14\p@\relax
  \bfseries
  \@title
  \end{center}%
} 
\makeatother
\title[Transience of a balanced excited random walk in $\ZZ^3$]{Transience of a balanced excited random walk\\ with vertical motion on return visits}
\date{}

\begin{document}
\begin{abstract}
We prove that the balanced excited random walk $M(3,2,1)$ is transient,
as conjectured by Benjamini, Kozma and Schapira. On its first departure
from each vertex of $\ZZ^3$, the walk takes a planar simple random walk
step; on every later departure, it takes a symmetric vertical step. We
also show that, almost surely, the range grows linearly in time up to
logarithmic factors. These results extend to walks whose first
departure is a planar simple random walk step with fixed probability
$\theta\in(0,1]$ and a symmetric vertical step otherwise.
\end{abstract}
\maketitle

\begin{figure}[!ht]
\centering
\makeatletter
\long\def\@makecaption#1#2{\vskip 5pt\relax
  {\small\raggedright\noindent{\scshape #1.}\enspace #2\par}}
\makeatother
\begin{minipage}[t]{.445\linewidth}
\centering
\begingroup
\definecolor{berwplanar}{HTML}{195970}
\definecolor{berwvertical}{HTML}{D9772E}
\definecolor{berwink}{HTML}{253E45}
\definecolor{berwgrid}{HTML}{C8D5D8}
\begin{tikzpicture}[x=1cm,y=1cm,line cap=round,line join=round,
  direction/.style={-{Stealth[length=4.6pt,width=3.8pt]},line width=1.20pt,
                    shorten >=2.3pt},
  guide/.style={draw=berwgrid,line width=.36pt},
  state/.style={font=\small\bfseries,text=berwink},
  site/.style={circle,draw=berwink,fill=white,inner sep=0pt,minimum size=3.6pt},
  idle/.style={circle,draw=berwink!45,fill=white,line width=.4pt,inner sep=0pt,
               minimum size=3.0pt},
  label/.style={font=\scriptsize,text=berwink}]
\path[use as bounding box] (0,0) rectangle (7.02,4.60);
\node[state] at (1.74,3.92) {First departure};
\node[state] at (5.27,3.92) {Later departures};
\foreach \xc/\mode in {1.74/0,5.27/1}{
 \begin{scope}[shift={(\xc,2.18)}]
 \begin{scope}[x={(1.06cm,-.41cm)},y={(-.81cm,-.53cm)},z={(0cm,1.02cm)}]
  \fill[berwgrid!19] (-.82,-.82,0)--(.82,-.82,0)--(.82,.82,0)--(-.82,.82,0)--cycle;
  \draw[guide] (-.82,-.82,0)--(.82,-.82,0)--(.82,.82,0)--(-.82,.82,0)--cycle;
  \draw[guide] (-1,0,0)--(1,0,0) (0,-1,0)--(0,1,0);
  \draw[guide,densely dashed] (0,0,-1)--(0,0,1);
  \ifnum\mode=0
    \foreach \tip in {{(1,0,0)},{(-1,0,0)},{(0,1,0)},{(0,-1,0)}}{
      \draw[direction,berwplanar] (0,0,0)--\tip;
      \node[site,draw=berwplanar,fill=berwplanar,minimum size=3.2pt] at \tip {};}
    \foreach \tip in {{(0,0,1)},{(0,0,-1)}} \node[idle] at \tip {};
    \node[label,anchor=south,yshift=2pt] at (0,0,0) {$v$};
  \else
    \foreach \tip in {{(0,0,1)},{(0,0,-1)}}{
      \draw[direction,berwvertical,line width=1.50pt] (0,0,0)--\tip;
      \node[site,draw=berwvertical,fill=berwvertical,minimum size=3.2pt] at \tip {};}
    \foreach \tip in {{(1,0,0)},{(-1,0,0)},{(0,1,0)},{(0,-1,0)}} \node[idle] at \tip {};
    \node[label,anchor=west,xshift=2pt] at (0,0,0) {$v$};
  \fi
  \node[site] at (0,0,0) {};
 \end{scope}
 \end{scope}
}
\node[label] at (1.74,.52) {$1/4$ per arrow};
\node[label] at (5.27,.52) {$1/2$ per arrow};
\end{tikzpicture}
\endgroup
\caption{At a vertex $v$, the first departure is planar and every
later departure is vertical. The filled neighbors in each panel are
the possible destinations, each with the indicated probability.}
\label{fig:model-rule}
\end{minipage}\hfill
\input{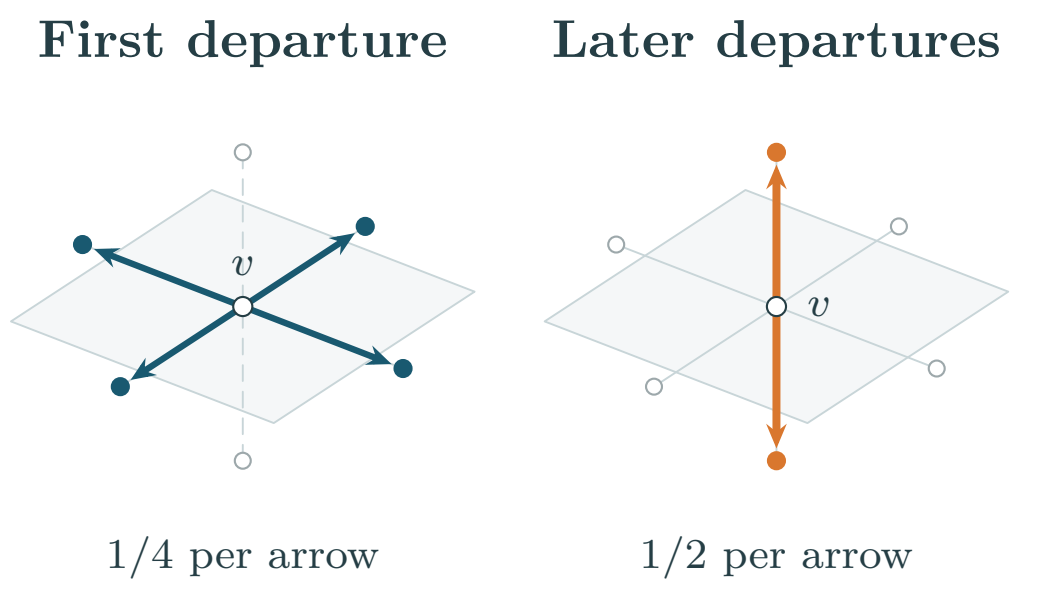}
\end{figure}

\section{Introduction}

\subsection{Background}

A walk on $\ZZ^d$ whose increments have conditional mean zero is a
martingale, but when the law of each step may depend on the past it need
not be transient, even for $d\ge3$. Peres, Popov and Sousi
\cite{PeresPopovSousi2013} constructed a recurrent walk in three
dimensions that chooses each step, by an adapted rule, from one of three
centered laws. They also proved that if every step is drawn from
one of two centered laws on $\RR^d$, $d\ge3$, with invertible covariance
matrices and finite moments of some order greater than two, then the walk
is transient, whatever the rule.

The balanced excited random walks of Benjamini, Kozma and Schapira
\cite{MR2788390} also choose between two centered laws, and each law may
be supported on a proper coordinate subspace. They were proposed as drift-free
relatives of the excited random walk of Benjamini and Wilson
\cite{MR1987097}, whose analysis relies on a bias in a fixed direction, and
as a step towards symmetric reinforced random walks. Let $e_1,\ldots,e_d$ be the standard basis of $\RR^d$, $d\ge2$, and
let $1\le d_1,d_2\le d$ be integers with $d_1+d_2\ge d$. The walk
$M(d,d_1,d_2)$ starts at the origin of $\ZZ^d$. On its first departure from
a vertex, it moves by a vector chosen uniformly from
$\{\pm e_i:1\le i\le d_1\}$. On every later departure from the same vertex,
it moves by a vector chosen uniformly from $\{\pm e_i:d-d_2+1\le i\le d\}$. Both laws
have mean zero, which explains the term \emph{balanced}. Benjamini, Kozma
and Schapira considered the \emph{non-overlapping} cases $d_1+d_2=d$. The
\emph{overlapping} cases $d_1+d_2>d$ were introduced by Camarena, Panizo and
Ram\'irez \cite{Camarena2020}, and our notation is adapted from theirs.

In two dimensions, Benjamini, Kozma and Schapira expressed the belief
that $M(2,1,1)$ is recurrent. Angel, Holmes and Ram\'irez \cite{MR4597323}
obtained upper and lower bounds for its range. The recurrence and range
of $M(2,1,2)$ are studied in \cite{Qin2026M212}. In a related planar model, each vertex
carries a horizontal or vertical label. Before each step the walker
switches the label at its current vertex with probability $q$, then takes
a symmetric step along the axis indicated by the updated label. Chan
\cite[Theorem 1.3]{Chan2023} proved that a single walker is almost surely
recurrent for $1/3<q\le1$ when the initial labels are independent and
uniform on $\{\mathrm H,\mathrm V\}$.

In the non-overlapping cases, the two coordinate projections, sampled
at their respective moves, follow simple random walks in dimensions
$d_1$ and $d_2$. Thus, the walk is transient whenever $\max\{d_1,d_2\}\ge3$, as pointed out by Benjamini, Kozma and Schapira. They also proved transience for $M(4,2,2)$. In dimensions at least four, Camarena, Panizo and
Ram\'irez proved transience in every overlapping case except $M(4,3,2)$. The second author studies the other nontrivial overlapping cases in dimension three, as well as $M(4,3,2)$, in \cite{Qin2026ERW3}.

In dimension three, Benjamini, Kozma and Schapira conjectured that both
$M(3,1,2)$ and $M(3,2,1)$ are transient
(see also \cite[Remark and Problem 3.22]{KosyginaZerner2013}).
For $M(3,1,2)$, whose first departures are one-dimensional and whose
subsequent departures are planar, Peres, Schapira and Sousi
\cite{MR3531697} proved the conjecture by a martingale defocusing argument.

\subsection{Main results}

In this paper we prove the conjecture for $M(3,2,1)$. We identify $\ZZ^3$
with $\ZZ^2\times\ZZ$, and call steps in the first two coordinates
\emph{planar} and steps in the third coordinate \emph{vertical}. Thus
$M(3,2,1)$ takes a planar simple random walk step on its first departure
from each vertex, and a symmetric vertical step on every later departure
(Figures \ref{fig:model-rule} and \ref{fig:model-simulation}).
Our proof also applies when the first departure is a fixed mixture of
planar and vertical steps. Fix $\theta\in[0,1]$, write $\delta_e$ for the
point mass at $e$, and let the laws of first and subsequent departures be
\begin{equation}\label{eq:fslaw}
 \mu_{\mathrm f}^{(\theta)}
 :=\frac\theta4\sum_{e\in\{\pm e_1,\pm e_2\}}\delta_e
   +\frac{1-\theta}{2}(\delta_{e_3}+\delta_{-e_3}),
 \qquad
 \mu_{\mathrm s}:=\frac12(\delta_{e_3}+\delta_{-e_3}).
\end{equation}
Write $S=(S_t)_{t\ge0}$ for the walk with these departure laws, started at
the origin with no previously visited vertices. The choice $\theta=1$ gives
$M(3,2,1)$. For $\theta=2/3$, the first departure is uniform over all six
nearest neighbors, so the walk is the overlapping model $M(3,3,1)$. For
$\theta=0$, the walk is a simple random walk on the vertical axis, and it
is recurrent. We call a walk \emph{transient} if, almost surely, every
vertex is visited only finitely often.

\begin{theorem}\label{mainthm}
For every $\theta\in(0,1]$, the walk with departure laws
\eqref{eq:fslaw} is transient. In fact, for every
$v\in\ZZ^3$,
\[
\EE\left[\sum_{t=0}^{\infty}\ind_{\{S_t=v\}}\right]<\infty.
\]
In particular, these conclusions hold for $M(3,2,1)$ and $M(3,3,1)$.
\end{theorem}

The proof also shows how the expected number of visits depends on the
distance from the starting point. Writing $\|x\|_1=|x_1|+|x_2|+|x_3|$,
we obtain, for every $a>0$ and every $x\in\ZZ^3$,
\begin{equation}\label{eq:green-decay}
 G(x):=\EE\!\left[\sum_{t\ge0}\ind_{\{S_t=x\}}\right]
 \le \frac{C_{\theta,a}}{[\log(2+\|x\|_1)]^a}.
\end{equation}
Our estimates also give a lower bound for the range. Write
\[
 \mathcal R_t:=\{S_s:0\le s<t\},\qquad R_t:=|\mathcal R_t|.
\]
Thus $\mathcal R_t$ consists of the vertices whose first departures have
occurred before time $t$. We call these vertices \emph{consumed}. A
vertex reached for the first time at time $t$ is \emph{fresh} and enters
$\mathcal R_{t+1}$ on its first departure.

\begin{proposition}\label{rangelowerbd}
For every $\theta\in(0,1]$ there is $c_\theta>0$ such that the walk
with departure laws \eqref{eq:fslaw} satisfies, almost surely,
\[
 R_t\ge c_\theta\frac{t}{(\log t)^{3/2}(\log\log t)^3}
\]
for all sufficiently large integers $t$.
\end{proposition}

For $M(d,d_1,d_2)$ we again write $R_t$ for the
number of distinct vertices visited before time $t$.
Since $R_t\le t$, Proposition \ref{rangelowerbd} gives
$\log R_t/\log t\to1$ almost surely for $M(3,2,1)$. For $M(3,1,2)$,
the stronger bound $R_t\ge c\,t/\log t$, eventually almost surely,
follows from the planar range law \cite[Theorem 4]{MR0047272}: either
$R_t\ge t/2$, or at least $\lfloor t/2\rfloor$ of the first $t$
departures are subsequent departures. The planar positions at these
departures form an initial segment of a planar simple random walk,
and all belong to the projection of $\mathcal R_t$.

However, the two models have small-range lower tails with different exponential orders. More generally, for $M(d,d_1,d_2)$, between two successive discoveries, once the first
departure has been made, the walk uses only the $d_2$ coordinates
available on subsequent departures. This suggests that $d_2$ should
govern the cost of remaining in a small set. The following Proposition \ref{prop-rangeLaplacelog} makes this precise. The proposition is independent of the transience argument and is proved directly from the departure rules in Section \ref{sec:range-lower-tail}.

 We write
$f(t)=\Theta(g(t))$ if $cg(t)\le f(t)\le Cg(t)$ for all large $t$, with
positive constants $c,C$ that depend only on the parameters indicated by
subscripts.

\begin{proposition}\label{prop-rangeLaplacelog}
For every fixed $\lambda>0$, the range $R_t$ of $M(d,d_1,d_2)$
satisfies, as $t\to\infty$,
\[
 -\log \EE[e^{-\lambda R_t}]
 =\Theta_{\lambda,d,d_1,d_2}
 \left(t^{d_2/(d_2+2)}\right).
\]
\end{proposition}

The exponent depends only on the number $d_2$ of coordinates used on
subsequent departures, including coordinates shared with the
first-departure law. It is $1/2$ for $M(3,1,2)$ and $1/3$ for $M(3,2,1)$.
For simple random walk on $\ZZ^d$, which is $M(d,d,d)$, it is the exponent
$d/(d+2)$ of Donsker and Varadhan \cite{DonskerVaradhan1979}. The
probability estimates \eqref{eq:range-tail} behind
Proposition \ref{prop-rangeLaplacelog} also compare the two
three-dimensional models at a common scale. They give
\begin{align*}
\PP\left(R_t\le t^{1/4}\right)
 &=\exp\left(-\Theta\!\left(t^{3/4}\right)\right)
 &&\text{for }M(3,1,2),\\
\PP\left(R_t\le t^{1/4}\right)
 &=\exp\left(-\Theta\!\left(t^{1/2}\right)\right)
 &&\text{for }M(3,2,1).
\end{align*}
Thus $M(3,2,1)$ is far more likely than $M(3,1,2)$ to remain in a small visited region for a long time.

\subsection{Why M(3,2,1) is harder than M(3,1,2)}
\label{sec:comparison}

In both models, a height martingale sampled along a planar walk
provides a route to transience. The difference lies in the control
available for its jumps and its bracket.

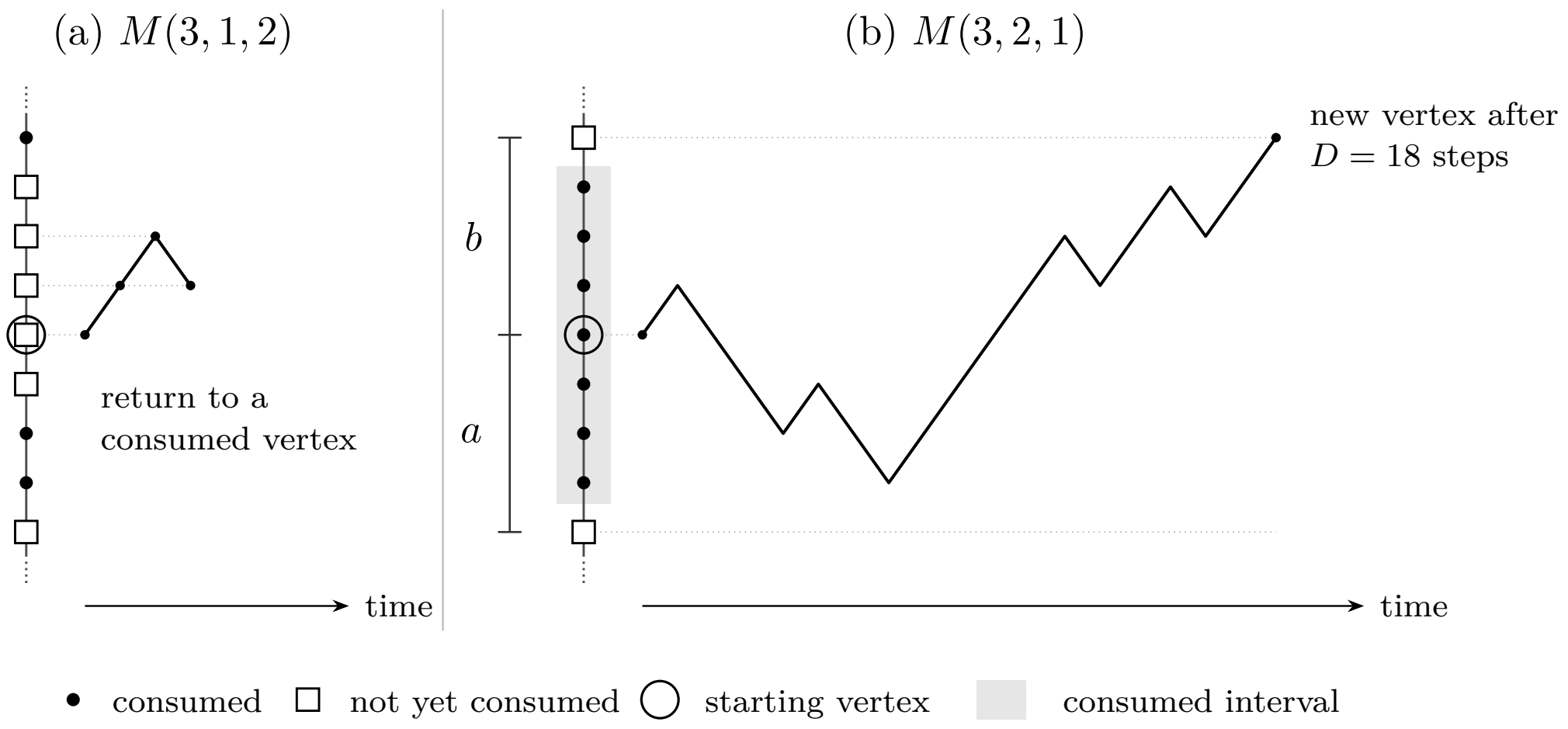
\begin{figure}[t]
\centering
\begin{tikzpicture}[
    x=1cm,y=0.42cm,
    >=Stealth,
    font=\small,
    consumed/.style={circle,fill=black,inner sep=0pt,minimum size=3.2pt},
    unvisited/.style={rectangle,draw=black,fill=white,line width=0.55pt,
                     inner sep=0pt,minimum size=5.4pt},
    walker/.style={circle,draw=black,line width=0.6pt,inner sep=0pt,
                   minimum size=9pt},
    run/.style={line width=0.75pt,line join=round},
    guide/.style={black!35,line width=0.3pt,densely dotted},
    axis/.style={line width=0.45pt,-{Stealth[length=4pt]}},
    lab/.style={font=\scriptsize}
]
\node at (1.85,6.1) {(a) $M(3,1,2)$};
\node at (8.6,6.1) {(b) $M(3,2,1)$};
\draw[black!25,line width=0.45pt] (4.15,-6.0) -- (4.15,6.6);

\def\xa{0.6}
\draw[black!70,line width=0.55pt] (\xa,-4.5) -- (\xa,4.5);
\draw[densely dotted,black!70,line width=0.55pt] (\xa,4.6) -- (\xa,5.1);
\draw[densely dotted,black!70,line width=0.55pt] (\xa,-4.6) -- (\xa,-5.1);
\foreach \h in {-3,-2,4}{\node[consumed] at (\xa,\h) {};}
\foreach \h in {-4,-1,0,1,2,3}{\node[unvisited] at (\xa,\h) {};}
\node[walker] at (\xa,0) {};
\draw[guide] (\xa+0.14,0) -- (1.1,0);
\draw[guide] (\xa+0.14,1) -- (2.0,1);
\draw[guide] (\xa+0.14,2) -- (1.7,2);
\draw[run] (1.1,0) -- (1.4,1) -- (1.7,2) -- (2.0,1);
\foreach \x/\h in {1.1/0,1.4/1,1.7/2,2.0/1}{\fill (\x,\h) circle[radius=1.15pt];}
\node[lab,anchor=north west,align=left] at (1.1,-0.75)
  {return to a\\consumed vertex};
\draw[axis] (1.1,-5.5) -- (3.35,-5.5) node[right,lab] {time};

\def\xb{5.35}
\fill[black!10] (\xb-0.23,-3.42) rectangle (\xb+0.23,3.42);
\draw[black!70,line width=0.55pt] (\xb,-4.5) -- (\xb,4.5);
\draw[densely dotted,black!70,line width=0.55pt] (\xb,4.6) -- (\xb,5.1);
\draw[densely dotted,black!70,line width=0.55pt] (\xb,-4.6) -- (\xb,-5.1);
\foreach \h in {-3,...,3}{\node[consumed] at (\xb,\h) {};}
\foreach \h in {-4,4}{\node[unvisited] at (\xb,\h) {};}
\node[walker] at (\xb,0) {};
\draw[black!80,line width=0.5pt] (4.72,-4) -- (4.72,4);
\foreach \h in {-4,0,4}{\draw[black!80,line width=0.5pt] (4.62,\h) -- (4.82,\h);}
\node[anchor=east] at (4.62,-2) {$a$};
\node[anchor=east] at (4.62,2) {$b$};
\foreach \h in {-4,4}{\draw[guide] (\xb+0.14,\h) -- (11.25,\h);}
\draw[guide] (\xb+0.14,0) -- (5.85,0);
\draw[run] plot coordinates {
  (5.85,0) (6.15,1) (6.45,0) (6.75,-1) (7.05,-2) (7.35,-1) (7.65,-2)
  (7.95,-3) (8.25,-2) (8.55,-1) (8.85,0) (9.15,1) (9.45,2) (9.75,1)
  (10.05,2) (10.35,3) (10.65,2) (10.95,3) (11.25,4)};
\fill (5.85,0) circle[radius=1.15pt];
\fill (11.25,4) circle[radius=1.15pt];
\node[lab,anchor=west,align=left] at (11.4,4)
  {new vertex after\\$D=18$ steps};
\draw[axis] (5.85,-5.5) -- (12.0,-5.5) node[right,lab] {time};

\node[consumed] at (1.0,-7.4) {};
\node[anchor=west,font=\scriptsize,text height=1.5ex,text depth=.25ex] at (1.2,-7.4) {consumed};
\node[unvisited] at (3.0,-7.4) {};
\node[anchor=west,font=\scriptsize,text height=1.5ex,text depth=.25ex] at (3.22,-7.4) {not yet consumed};
\node[walker] at (6.0,-7.4) {};
\node[anchor=west,font=\scriptsize,text height=1.5ex,text depth=.25ex] at (6.25,-7.4) {starting vertex};
\fill[black!10] (8.7,-7.8) rectangle (9.12,-7.0);
\node[anchor=west,font=\scriptsize,text height=1.5ex,text depth=.25ex] at (9.3,-7.4) {consumed interval};
\end{tikzpicture}
\caption{The one-dimensional coordinate in each model, plotted against time on
the same scales. The markers at the left of each panel show the initial
configuration. In (a) the third step reverses the second, returning to a
vertex already visited during the run. In (b) the nearest unvisited heights
are $a=b=4$ steps from the starting point, and the walk exits the shaded
interval after $D=18$ vertical steps.}

\label{fig:vertical-motion}
\end{figure}

For $M(3,1,2)$, call the first coordinate the height and the last two
coordinates the planar position. Peres, Schapira and Sousi
\cite{MR3531697} observe
the height after successive planar moves. Given the planar path, the
sampled heights form a martingale. During a run of one-dimensional
steps, a reversal takes the walk back to a vertex it has just left,
where the next step must be planar, see Figure \ref{fig:vertical-motion} (a). The run length therefore has a geometric tail uniformly in the past. This gives the geometric bound on the height jumps required by their defocusing theorem and bounds
the increase of the bracket (the sum of the conditional variances of
the height increments) by $Ck$ over any $k$ planar moves. For the
lower bound, the conditional variance of each height increment equals
the conditional mean number of one-dimensional steps before the next
planar move. Every newly visited planar site forces at least one such
step and so contributes at least one unit to the bracket. Conditional
on the planar walk returning to the origin, their planar range
estimates supply the bracket lower bounds needed for defocusing.

For $M(3,2,1)$ the height is the third coordinate, and we record it
at successive discoveries of new vertices. After a planar step, the
walk stays in its column (the vertical line through the point of arrival)
until it discovers a new vertex. Call this vertical segment an
\emph{excursion} and denote its number of steps by $D$. If the planar
step itself discovers a vertex, then $D=0$. Otherwise the walk may
reverse direction many times before it leaves the interval of consumed
heights containing the point of arrival, see Figure \ref{fig:vertical-motion} (b).
Given the arrival configuration, let $a$ and $b$ be the distances to
the lower and upper exits of this interval. The conditional variance of
the height change and the conditional mean of $D$ are both $ab$. 
The departure rule thus gives neither a uniform geometric tail for
the height jumps nor a uniform bound on the conditional mean excursion
duration. The defocusing theorem of Peres, Schapira and Sousi does not
apply directly, as they observed in \cite[Section 4]{MR3531697}.
Our Proposition \ref{lem:rectangle} treats the jump bound separately
from the bounds on the bracket. The main work is to control the
vertical time accumulated over intervals of discovery indices.

\subsection{Outline of the proof}

We describe the proof for transience of $M(3,2,1)$. Write $(U_j,H_j)$ for the planar
position and height at the $j$th discovery, counting the origin as
discovery zero. From one discovery to the next there is exactly one
planar step, so $U$ is a planar simple random walk. For fixed
$u\in\ZZ^2$, the probability of $U_n=u$ is of order $1/n$ for those $n$
with $\PP(U_n=u)>0$. To obtain a summable bound on visits to a
vertex $(u,z)$, we therefore estimate the spread of $H_n$ under the
planar bridge law $Q_n^u=\PP(\,\cdot\mid U_n=u)$.

Let $T_j$ count the vertical steps up to the $j$th discovery and fix
$q>1$. For an excursion starting in an interval of length $k$ consisting of consumed vertices,
the conditional $L^q$ norm of its duration, given the starting
configuration, is at most $C_q(k+1)^2$. The interval length is at most
the number of visits by $U$ to that column before index $n$. At any
fixed excursion index at most $n$, this visit count has
$L^{2q}(Q_n^u)$ norm at most $C_{q,u}\log(n+2)$, see Lemmas \ref{lem:height-martingale} and \ref{lem:bridge-local}. By Minkowski's inequality, summing over the last $r$ excursions
therefore gives, for $1\le r\le n$,
\begin{equation}\label{eq:crude-clock}
 \|T_n-T_{n-r}\|_{L^q(Q_n^u)}\le C_{q,u}\,r[\log(n+2)]^2 .
\end{equation}
The square of the logarithm in \eqref{eq:crude-clock} is too large
to yield a summable bound by Proposition \ref{lem:rectangle}. To reduce the power of $\log n$, we count together the vertical steps in each column.

For a deterministic interval in a fixed column, the mean exit time
of a one-dimensional simple random walk, viewed as a function of the
starting height, decreases by one on average at each vertical step
inside the interval.
Compensating for planar steps into and out of the column by a weighted
count of first departures from the column and its four neighbors gives
the occupation identity of
Lemma \ref{lem:column-identity}. At discovery $n$, however, the maximal
intervals of consumed heights have random endpoints. Using these
endpoints as coefficients of earlier martingale increments would make
the coefficients depend on the future. Fourier inversion separates
the interval endpoints from random sums whose coefficients depend
only on the history at the step in question. Under the bridge law,
exchangeability of the remaining planar increments controls the bias
arising from first departures. Lemma \ref{lem:column-occupation} then
gives an $L^q$ bound of order $[\log(n+2)]^{5/2}$ for the total vertical
time spent in a fixed column up to discovery $n$.

To sum over the columns visited in the last $r$ excursions, we must
account for how the planar walk selects them. We split according to
the last visit to each column in the window. For a fixed index $j$ to
be the last visit, the later planar increments must avoid a return;
on this event the total number of vertical steps in that column up to discovery $n$ is already determined at discovery $j$. Summing over the possible last visits gives the
factor $r/\log(r+2)$ from the expected range of a planar walk of $r$
steps. Consequently, on windows of length $r\ge n^\alpha$ with fixed
$0<\alpha<1$, Proposition \ref{prop:window} gives
\[
 \|T_n-T_{n-r}\|_{L^q(Q_n^u)}
 \le C_{q,u}\,\frac{r[\log(n+2)]^{5/2}}{\log(r+2)}
 \le C_{q,u,\alpha}\,r[\log(n+2)]^{3/2}.
\]

To apply the defocusing estimate, reveal the entire planar path.
The unused vertical steps still have independent symmetric signs,
so the heights at discovery times form a martingale with this
additional information. Each bracket increment is the conditional
mean duration of the corresponding vertical excursion. For the last
$r$ excursions, inequality \eqref{eq:predictable-projection} bounds
the $L^q(Q_n^u)$ norm of the bracket increase by
$q\|T_n-T_{n-r}\|_{L^q(Q_n^u)}$. Thus the preceding estimate for
$T_n-T_{n-r}$ provides the required upper bounds on sums of conditional
variances. Starting halfway through the discovery indices, we form
successive blocks by halving the number of indices remaining until
$n$, stopping when this number is of order $n^{1/2}$. The upper
bracket bounds are used on the longer intervals from the beginning
of each block to discovery $n$. An immediate reversal of a planar
step forces a vertical step in at least one of the two corresponding
excursions. With high probability under $Q_n^u$, enough such pairs
occur in each block to give a bracket lower bound proportional to
its length. The height jumps are bounded by the maximum number of
visits of $U$ to a site before index $n$. Applying martingale
defocusing over the order $\log n$ blocks, and adding the probabilities
that the bracket and jump conditions fail, gives, for every fixed
$A>0$ (Proposition \ref{prop:height-smallball}),
\[
 \sup_{z\in\RR} Q_n^u\bigl(|H_n-z|\le n^{1/4}\bigr)
 \le C_{A,u}[\log(n+2)]^{-A}.
\]

We next turn this height estimate into a bound on visits at all
physical times. Conditional on $U_n=u$, excursion $n$ can visit
$(u,z)$ only if $|H_n-z|\le n^{1/4}$ or it exits an interval of
length greater than $n^{1/4}$. The interval length is bounded by
the number of visits of $U$ to $u$ before index $n$. The tail estimate
for this visit count shows that the probability of the second event
also decays faster than any inverse power of $\log n$.
The same decay therefore holds for the conditional probability that
excursion $n$ visits $(u,z)$. The number of visits during the excursion
is at most its duration plus one, whose $L^2(Q_n^u)$ norm is
$O((\log n)^2)$. Choosing $A$ sufficiently large, Cauchy--Schwarz
and the bound $\PP(U_n=u)=O(1/n)$ give an upper bound of
$O(1/(n(\log n)^2))$ for the expected number of visits contributed
by excursion $n$. Summing over $n$ proves transience.

The argument extends to $0<\theta<1$, with constants depending
on $\theta$.

Section \ref{secpreli} constructs the walk from independent
planar and vertical walks and proves the required bridge estimates.
Section \ref{sec:column} gives the column occupation identity and
the terminal-window estimate. Section \ref{sec:defocusing} proves
the martingale defocusing estimate. These results are combined in
Section \ref{sec:transience} to prove transience, finite expected
occupation times, and the range results.

\section{Notation and preliminary estimates}
\label{secpreli}

\subsection{The process and basic estimates}

Fix $\theta\in(0,1]$ and write a point of $\ZZ^3=\ZZ^2\times\ZZ$ as a
planar position and a height. We use $t$ for the time of the original
walk and $j,m,n$ to count discoveries. All these indices are integers.
Throughout the paper, $c,C$ denote positive constants that may change
from line to line and may depend on $\theta$. Subscripts indicate
further dependence. We write $L_n=\log(n+2)$ and often set $L=L_n$.

We construct $S$ from two independent walks. The planar walk
$(U_j)_{j\ge0}$ starts at $0$ and has independent increments of law
\begin{equation}\label{eq:independent-drivers}
 \nu_\theta(0):=1-\theta,\qquad
 \nu_\theta(e):=\frac\theta4,\quad
 e\in\{(\pm1,0),(0,\pm1)\}.
\end{equation}
Set $\nu_\theta(e)=0$ for all other $e\in\ZZ^2$.
Let $(\widetilde H_k)_{k\ge0}$ be an independent one-dimensional
simple random walk started at $0$, with signs
$\tilde\epsilon_k=\widetilde H_k-\widetilde H_{k-1}$ for $k\ge1$.
Write $\PP$ for the joint law of $U$ and $\widetilde H$, and hence of
the walk $S$ constructed from them.

Each time $S$ leaves a vertex for the first time, read the next increment
of $U$. For $\theta=1$ this always gives a planar step. For $0<\theta<1$,
a nonzero increment gives the indicated planar step, while a zero increment
calls for a vertical step using the next increment of $\widetilde H$. At every
subsequent departure from that vertex, read only the next vertical
increment. Each choice uses only information already revealed, so the
unused increments give precisely the transition probabilities
in \eqref{eq:fslaw}. In particular, the number of planar increments
read by time $t$ is the number $R_t$ of first departures before $t$,
including those that produced a vertical step. Thus the planar position
of $S_t$ is $U_{R_t}$. For $\theta=1$, the gray trace in
Figure \ref{fig:model-simulation} shows this projection.

Let $\FF_t:=\sigma(S_0,\ldots,S_t)$. A vertex is \emph{consumed} at time
$t$ if its first departure has already occurred, that is, if it belongs
to $\mathcal R_t=\{S_0,\ldots,S_{t-1}\}$, with $\mathcal R_0=\varnothing$.
The current vertex is \emph{fresh} if
$S_t\notin\mathcal R_t$. The discovery times are
\[
 \tau_0=0,\qquad
 \tau_{j+1}=\inf\{t>\tau_j:S_t\notin\mathcal R_t\}.
\]
The first-departure steps are precisely those from $\tau_j$ to $\tau_j+1$.
Each $\tau_j$ is finite almost surely. Indeed, after a first departure,
either the arrival is fresh or the subsequent vertical walk exits a finite
interval of consumed vertices. At time $\tau_j$, exactly $j$ vertices were
visited earlier, so $R_{\tau_j}=j$.

Write $W_j=S_{\tau_j}=(U_j,H_j)$, where $H_j$ is the height at discovery
$j$ and $W_0=0$. Let $\kappa_t$ count the vertical steps made before
physical time $t$, and let $T_j$ be this count at discovery $j$. Thus
\[
 \kappa_t=\sum_{0\le r<t}\ind_{\{S_{r+1}-S_r\in\{e_3,-e_3\}\}}
 \qquad\text{and}\qquad T_j=\kappa_{\tau_j}.
\]
The step arriving at $W_j$ is included if it is vertical. In particular
$H_j=\widetilde H_{T_j}$, with a clock $T_j$ that depends on both walks.
For $j\ge1$, let $D_j$ count the vertical steps between
$\tau_{j-1}$ and $\tau_j$, including a possible first vertical departure. Then
\begin{equation}\label{eq:clock}
 \begin{split}
 D_j&=T_j-T_{j-1}
     =\tau_j-\tau_{j-1}-\ind_{\{U_j\ne U_{j-1}\}},\qquad
 T_n=\sum_{j=1}^nD_j,\\
 \tau_n&=T_n+\sum_{j=1}^n\ind_{\{U_j\ne U_{j-1}\}}
 \le T_n+n.
 \end{split}
\end{equation}

\begin{figure}[t]
\centering
\begin{tikzpicture}[
    x=1cm,y=0.67cm,
    >=Stealth,
    font=\small,
    consumed/.style={circle,fill=black,inner sep=0pt,minimum size=3.2pt},
    unvisited/.style={rectangle,draw=black,fill=white,line width=0.55pt,
                     inner sep=0pt,minimum size=5.4pt},
    first/.style={->,line width=0.7pt,
                  >={Stealth[length=1.7mm,width=1.15mm]}}
]

\draw[black!25,line width=0.45pt] (7.65,-3.9) -- (7.65,4.75);
\node at (3.55,4.48) {(a) $U_j\ne U_{j-1}$};
\node at (11.45,4.48) {(b) $U_j=U_{j-1}$};

\def\xleft{3.0}
\fill[black!10] (2.77,-2.24) rectangle (3.23,2.24);
\draw[black!70,line width=0.55pt] (\xleft,-3.48) -- (\xleft,3.48);
\foreach \h in {-2,-1,0,1,2}{\node[consumed] at (\xleft,\h) {};}
\foreach \h in {-3,3}{\node[unvisited] at (\xleft,\h) {};}
\draw[densely dotted,black!70,line width=0.55pt] (\xleft,3.52) -- (\xleft,3.98);
\draw[densely dotted,black!70,line width=0.55pt] (\xleft,-3.52) -- (\xleft,-3.98);
\node[anchor=east] at (2.80,3.0) {$W_j$};
\node at (\xleft,-4.48) {$U_j$};

\node[consumed] at (0.45,0) {};
\node[anchor=north west] at (0.30,-0.28) {$W_{j-1}$};
\draw[first] (0.58,0) -- (2.90,0);
\node[font=\scriptsize,anchor=south,align=center] at (1.61,0.30) {first planar\\departure};
\node[anchor=north east] at (2.84,-0.24) {$H_{j-1}$};

\draw[black!80,line width=0.5pt] (3.85,-3) -- (3.85,3);
\foreach \h in {-3,0,3}{\draw[black!80,line width=0.5pt] (3.73,\h)--(3.97,\h);}
\node[anchor=west] at (4.05,-1.5) {$a$};
\node[anchor=west] at (4.05,1.5) {$b$};
\draw[decorate,decoration={brace,amplitude=3.5pt}] (4.75,3) -- (4.75,0)
  node[midway,right=6pt] {$H_j-H_{j-1}=b$};

\def\xright{10.85}
\fill[black!10] (10.62,-2.24) rectangle (11.08,2.24);
\draw[black!70,line width=0.55pt] (\xright,-3.48) -- (\xright,3.48);
\foreach \h in {-2,-1,0,1,2}{\node[consumed] at (\xright,\h) {};}
\foreach \h in {-3,3}{\node[unvisited] at (\xright,\h) {};}
\draw[densely dotted,black!70,line width=0.55pt] (\xright,3.52) -- (\xright,3.98);
\draw[densely dotted,black!70,line width=0.55pt] (\xright,-3.52) -- (\xright,-3.98);
\node[anchor=east] at (10.65,3.0) {$W_j$};
\node at (\xright,-4.48) {$U_j=U_{j-1}$};
\node[anchor=east] at (10.48,-1.0) {$W_{j-1}$};

\draw[first,line width=0.85pt] (10.92,-0.93) to[bend right=26] (10.92,-0.07);
\node[font=\scriptsize,anchor=east,align=right] at (10.46,-0.20)
  {first vertical\\departure};

\draw[black!80,line width=0.5pt] (11.70,-3) -- (11.70,3);
\foreach \h in {-3,-1,3}{\draw[black!80,line width=0.5pt] (11.58,\h)--(11.82,\h);}
\node[anchor=west] at (11.90,-2) {$a$};
\node[anchor=west] at (11.90,1) {$b$};
\draw[decorate,decoration={brace,amplitude=3.5pt}] (12.60,3) -- (12.60,-1)
  node[midway,right=6pt] {$H_j-H_{j-1}=b$};

\node[consumed] at (2.15,-5.18) {};
\node[anchor=west,font=\scriptsize] at (2.38,-5.18) {consumed};
\node[unvisited] at (5.72,-5.18) {};
\node[anchor=west,font=\scriptsize] at (5.97,-5.18) {not yet consumed};
\fill[black!10] (9.87,-5.36) rectangle (10.29,-5.00);
\node[anchor=west,font=\scriptsize] at (10.49,-5.18) {consumed interval};
\end{tikzpicture}
\caption{The $j$th vertical excursion, shown exiting above the consumed
interval. The shaded interval is taken in column $U_j$ after the first
departure from $W_{j-1}$. The distances $a$ and $b$ are measured from
$H_{j-1}$ to the sites immediately below and above the shaded interval,
respectively. In (b) the vertical first departure consumes $W_{j-1}$ and counts
toward $D_j$. In the configuration shown, adding $W_{j-1}$ joins two
previously separate intervals.
Only the first departure is drawn.}
\label{fig:column-excursion}
\end{figure}
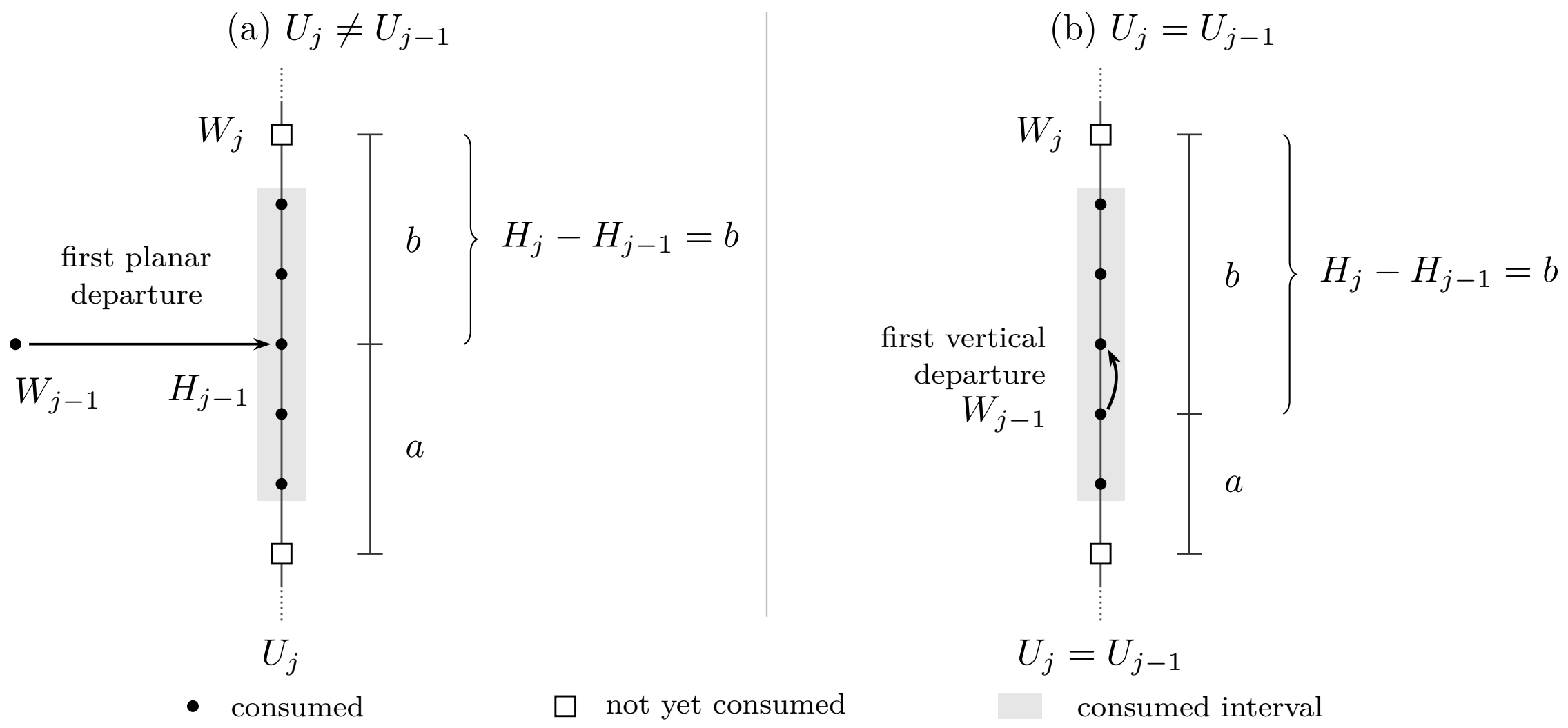

For $j\ge1$, we call the vertical segment between the discoveries of
$W_{j-1}$ and $W_j$ the \emph{$j$th excursion},
see Figure \ref{fig:column-excursion} for an illustration. If $U_j\ne U_{j-1}$, it starts
at the planar arrival $(U_j,H_{j-1})$ and occupies the time interval
$[\tau_{j-1}+1,\tau_j]$. If $U_j=U_{j-1}$, it starts at $W_{j-1}$
and occupies $[\tau_{j-1},\tau_j]$, including the first vertical step.
In either case it stays in column $U_j$ and ends at $W_j$. Counting
both endpoints gives $D_j+1$ integer times, including the single time
$\tau_j$ when a planar step discovers a fresh vertex and $D_j=0$.

At the discovery time $\tau_n$, the consumed set and its column counts are
\begin{equation}
 \begin{split}
 \mathcal R_{\tau_n}&=\{W_0,\ldots,W_{n-1}\},\qquad R_{\tau_n}=n,\\
 \ell_n(v)&:=\sum_{j=0}^{n-1}\ind_{\{U_j=v\}}
       =|\mathcal R_{\tau_n}\cap(\{v\}\times\ZZ)|,\qquad
 \ell_n^*:=\max_{v\in\ZZ^2}\ell_n(v).
 \end{split}
\end{equation}
Thus $\ell_n(v)$ is the local time of $U$ at $v$ before index $n$, and
$\ell_n^*$ is its maximum over $v$. The current vertex $W_n$ has not yet been
consumed and is not included in these counts.

Write $Y_j=U_j-U_{j-1}$ for the planar increment read at the $j$th
first departure. At discovery $j$ we know the history
$\mathcal D_j:=\FF_{\tau_j}$. Just before the vertical part of excursion
$j$, the information we use is
\[
 \mathcal A_j:=\mathcal D_{j-1}\vee\sigma(Y_j).
\]
The past determines the starting height $H_{j-1}$, and $Y_j$ selects
the column $U_j$. When $Y_j=0$, the first departure is itself vertical.
Conditioning on $\mathcal A_j$ leaves its sign unrevealed.

To describe the excursion, adjoin $W_{j-1}$ to the set consumed before
its first departure. In column $U_j$, let $I_j$ be the maximal interval of
consecutive heights in this enlarged set containing $H_{j-1}$. Take
$I_j=\varnothing$ if $(U_j,H_{j-1})$ is outside the enlarged set.
Thus $I_j$ is $\mathcal A_j$-measurable and $|I_j|\le\ell_j(U_j)$.
It is empty precisely when a planar step discovers a fresh vertex.

\begin{lemma}\label{lem:height-martingale}
Let $Q=\PP$ or $Q=\PP(\,\cdot\mid U_n=u)$ with $\PP(U_n=u)>0$.
\begin{enumerate}
\item[(a)] If $\rho$ is a $Q$-almost surely finite stopping time for
$(\FF_t)$, then, conditionally on $\FF_\rho\vee\sigma(U)$, the unused
vertical signs $(\tilde\epsilon_{\kappa_\rho+i})_{i\ge1}$ are independent
and take each of the values $-1$ and $1$ with probability $1/2$ under $Q$.
\item[(b)] Given $\mathcal A_j$, the height in excursion $j$ is a simple
random walk started at $H_{j-1}$ and stopped on leaving $I_j$. If
$I_j=\varnothing$, then $D_j=0$ and $H_j=H_{j-1}$. Otherwise, let $a$
and $b$ be the distances from $H_{j-1}$ to the first heights outside
$I_j$ below and above it. Then
\begin{equation}\label{eq:exit-identities}
 \EE_Q[H_j-H_{j-1}\mid\mathcal A_j]=0
 \quad\text{and}\quad
 \EE_Q[(H_j-H_{j-1})^2\mid\mathcal A_j]
     =\EE_Q[D_j\mid\mathcal A_j]=ab .
\end{equation}
For every integer $k\ge0$, on $\{|I_j|\le k\}$ we have $|H_j-H_{j-1}|\le k$ and
\begin{equation}\label{eq:crude-exit}
 \bigl(\EE_Q[D_j^p\mid\mathcal A_j]\bigr)^{1/p}
 \le C p(k+1)^2\qquad(p\ge1).
\end{equation}
\item[(c)] The heights $(H_j)_{j\ge0}$ form a square-integrable
martingale under $Q$ in the discovery filtration $(\mathcal D_j)$, with
\begin{equation}
 \EE_Q[(H_{j+1}-H_j)^2\mid\mathcal D_j]
 =\EE_Q[D_{j+1}\mid\mathcal D_j]\qquad(j\ge0).
\end{equation}
The conclusions in (b) and (c) also hold with $\sigma(U)$ adjoined to
$\mathcal A_j$ and $\mathcal D_j$.
\end{enumerate}
\end{lemma}

\begin{proof}
For (a), the planar endpoint condition involves only $U$, so under
either law the vertical signs are independent symmetric signs,
independent of $U$. Let
$\mathcal K_k=\sigma(U,\tilde\epsilon_1,\ldots,\tilde\epsilon_k)$.
Given $U$ and the first $k$ signs, run the walk until time $\rho$ is
reached or the next sign is needed. This determines $\{\kappa_\rho\le k\}$
and, on $\{\kappa_\rho=k\}$, the history through $\rho$. Thus $\kappa_\rho$
is a stopping time for $(\mathcal K_k)$ and
$\FF_\rho\vee\sigma(U)\subseteq\mathcal K_{\kappa_\rho}$.
For $A\in\mathcal K_{\kappa_\rho}$, the event
$A\cap\{\kappa_\rho=k\}$ belongs to $\mathcal K_k$ and is independent
of the signs after $k$. Summing over $k$ proves that the signs after
$\kappa_\rho$ have their original product law given $\mathcal K_{\kappa_\rho}$.

For (b), apply (a) at $\rho=\tau_{j-1}$. Given
$\mathcal D_{j-1}\vee\sigma(U)$, the starting height and $I_j$ are
known. If $I_j$ is nonempty, the walk reads consecutive vertical signs
until it reaches one of the two adjacent heights outside $I_j$. Both
are fresh in this column, so this first exit ends the excursion. When
$Y_j=0$ the exit time includes the first vertical departure. This
description depends only on $H_{j-1}$ and $I_j$, hence it holds also
given the smaller sigma-field $\mathcal A_j$.

Write $X_s$ for the displacement after $s$ vertical steps of the simple
walk in $I_j$. A run of $|I_j|$ positive signs forces an exit, so $D_j$
is finite. Stop $X_s$ and $X_s^2-s$ at $D_j\wedge m$. The stopped
displacements stay in $[-a,b]$, and $D_j\wedge m$ increases to $D_j$.
Letting $m\to\infty$ gives $\EE_Q[X_{D_j}\mid\mathcal A_j]=0$ and
$\EE_Q[D_j\mid\mathcal A_j]=\EE_Q[X_{D_j}^2\mid\mathcal A_j]$.
The exit displacement equals $-a$ or $b$, with probabilities $b/(a+b)$
and $a/(a+b)$, so its second moment is $ab$. This proves
\eqref{eq:exit-identities}.

If $1\le |I_j|\le k$, then $a+b=|I_j|+1$, so the exit displacement has
absolute value at most $k$ and the mean exit time from every height
in $I_j$ is at most $(k+1)^2/4$. Markov's inequality and the Markov
property at multiples of $(k+1)^2$ give
\[
 Q\bigl(D_j>m(k+1)^2\mid\mathcal A_j\bigr)\le4^{-m}
 \qquad(m=0,1,\ldots).
\]
Integrating this tail gives \eqref{eq:crude-exit}. 

Finally, $|I_j|\le j$ gives square integrability at each discovery.
Conditioning \eqref{eq:exit-identities} on $\mathcal D_{j-1}$ proves
(c), with both the duration and the displacement zero on an empty
interval. The same exit law was proved given
$\mathcal D_{j-1}\vee\sigma(U)=\mathcal A_j\vee\sigma(U)$, which gives
the assertions with the full planar path revealed.
\end{proof}

\subsection{Estimates for the planar walk}\label{sec:bridge}

The positions $W_0,\ldots,W_j$ reveal $U_0,\ldots,U_j$. Conversely, the first $j$
planar increments and the entire path of $\widetilde H$ suffice to run the construction
through $\tau_j$. Thus
\[
 \sigma(Y_1,\ldots,Y_j)\ \subseteq\ \mathcal D_j\
 \subseteq\ \sigma(Y_1,\ldots,Y_j,\widetilde H).
\]
This will allow us to condition on the observed three-dimensional walk
while using the independence of the planar increments still to come
under $\PP$.

Let $p_n(u)=\PP(U_n=u)$ and call $(n,u)$ \emph{admissible} when
$p_n(u)>0$. The planar local central limit theorem
\cite[Theorems 2.1.1 and 2.1.3]{MR2677157} gives
\begin{equation}\label{eq:planar-kernel-uniform}
 \sup_y p_n(y)\le \frac C{n+1},\qquad
 p_n(u)\ge\frac{c_u}{n+1}\quad\text{for admissible }(n,u).
\end{equation}
The second bound is for each fixed $u$, with $c_u$ decreased if necessary
to include small admissible times. When $\theta=1$ the walk has period
two and the periodic form of the theorem applies. Summing the lower bound
at $u=0$ gives $\sum_{i=0}^n p_i(0)\ge c L_n$, using the even indices
when $\theta=1$.

We will estimate joint events involving $U_n=u$ and later recover
conditional estimates when $u$ is fixed. For any event $E$ define
\[
 \Pi_n^u(E)=(n+1)\PP(E\cap\{U_n=u\}),\qquad
 w_{n,u}=\Pi_n^u(\Omega)=(n+1)p_n(u).
\]
At an admissible endpoint this finite measure is $w_{n,u}Q_n^u$, where
$Q_n^u=\PP(\,\cdot\mid U_n=u)$ is the bridge law. By
\eqref{eq:planar-kernel-uniform}, $w_{n,u}\le C$ for every $n,u$ and
$w_{n,u}\ge c_u$ at the admissible times for each fixed $u$. Thus a
bound under $\Pi_n^u$ can be uniform in $u$, and
$\|X\|_{p,Q_n^u}=w_{n,u}^{-1/p}\|X\|_{p,\Pi_n^u}$ recovers the
corresponding bound for a fixed endpoint. When $p_n(u)=0$, $\Pi_n^u$
is the zero measure and any upper bound for its integrals holds trivially.

In this subsection and Section \ref{sec:column} we use
$(Q,\mu)=(\PP,\PP)$ for the unconditioned walk and
$(Q,\mu)=(Q_n^u,\Pi_n^u)$ for an admissible bridge. Expectations
conditional on a history are taken under the probability law $Q$.
For event masses and moments we use $\mu$ and write
$\|X\|_{p,\mu}=(\int |X|^p\,d\mu)^{1/p}$.

To compare the first half of a bridge with the unconditioned walk, let
$0\le a\le\lceil n/2\rceil$ and let $Z$ be a nonnegative function of
$Y_1,\ldots,Y_a$ and $\widetilde H$. Under $\PP$, the remaining $n-a$ planar increments are independent of
$Y_1,\ldots,Y_a$ and $\widetilde H$. Given the first $a$ increments
and the vertical walk, the conditional probability of $U_n=u$ is
therefore $p_{n-a}(u-U_a)$. This gives
\begin{equation}\label{eq:initial-density}
 \int Z\,d\Pi_n^u
 =\EE\big[Z(n+1)p_{n-a}(u-U_a)\big]\le C\EE Z.
\end{equation}
Here \eqref{eq:planar-kernel-uniform} bounds the density by
$C(n+1)/(n-a+1)\le 2C$. The comparison applies in particular to
$\mathcal D_a$-measurable variables. To read the path from the other
endpoint, put $\widehat U_i=U_n-U_{n-i}$ for $0\le i\le n$. Its
increments are $Y_n,\ldots,Y_1$. Permuting the increments preserves
their product probability and the event that their sum is $u$. Thus
$(\widehat U_i)_{0\le i\le n}$ and $(U_i)_{0\le i\le n}$ have the
same law under the bridge, and for every
nonnegative function $\Phi$,
\begin{equation}\label{eq:block-density}
 \int \Phi(Y_n,\ldots,Y_{n-a+1})\,d\Pi_n^u
       \le C\,\EE\Phi(Y_1,\ldots,Y_a).
\end{equation}

\begin{lemma}\label{lem:bridge-local}
For either choice of $\mu$, every $x\ge0$, $v\in\ZZ^2$ and $p\ge1$,
\begin{align}
 \mu(\ell_n(v)>x)&\le C e^{-cx/L_n},&
 \mu(\ell_n^*>x)&\le C(n+1)^2e^{-cx/L_n},\label{eq:local-tails}\\
 \sup_{0\le j\le n}\|\ell_n(U_j)+1\|_{p,\mu}&\le Cp L_n,&
 \|\ell_n^*\|_{p,\mu}&\le C L_n(L_n+p).\label{eq:local-index}
\end{align}
In particular $\sup_v\|\ell_n(v)+1\|_{p,\mu}\le Cp L_n$.
\end{lemma}

\begin{proof}
We first give a counting argument that will also apply to the planar
range. Let $I_0,\ldots,I_m$ be indicator variables adapted to a
filtration $(\mathcal G_i)$ and let $N=\sum_{i=0}^m I_i$. Suppose that
$A\ge1$, $\EE N\le A$, and
$\EE[\sum_{l=i+1}^m I_l\mid\mathcal G_i]\le A$ at each index $i$.
The binomial coefficient $\binom Nk$ counts the increasing $k$-tuples
of indices at which the indicator is one. Summing over the first $k-1$
indices and conditioning at the last of them gives, for $k\ge2$,
\[
 \EE\binom Nk
 =\sum_{i_1<\cdots<i_{k-1}}
   \EE\!\left[I_{i_1}\cdots I_{i_{k-1}}\,
     \EE\!\left[\sum_{l=i_{k-1}+1}^m I_l\,
                    \middle|\,\mathcal G_{i_{k-1}}\right]\right]
 \le A\,\EE\binom N{k-1}.
\]
Induction gives $\EE\binom Nk\le A^k$. The binomial expansion now yields
\[
 \EE e^{N/(4A)}\le\EE(1+1/(2A))^N
       =\sum_{k\ge0}(2A)^{-k}\EE\binom Nk\le2,
\]
where we used $\log(1+1/(2A))\ge 1/(4A)$. Markov's inequality gives
$\PP(N>x)\le2e^{-x/(4A)}$. Also
$y^p\le (p/e)^p e^y$ for $y\ge0$, so $\|N\|_{p,\PP}\le CpA$ for
every real $p\ge1$.

For visits to a fixed site $v$, take $I_i=\ind_{\{U_i=v\}}$ and the
natural filtration of $U$. Conditional on the walk through index $i$,
the expected number of visits between $i$ and $m$, including both ends, is
\[
 \sum_{l=0}^{m-i}p_l(v-U_i)
       \le C\sum_{l=0}^m\frac1{l+1}\le C\log(m+2).
\]
The same bound holds for the mean number of visits from any starting
point. The counting argument proves the fixed-site tail under $\PP$,
including visits at the last index if desired.

For a bridge, set $a=\lfloor n/2\rfloor$. Since
$U_{n-i}=u-\widehat U_i$ on $\{U_n=u\}$,
\[
 \ell_n(v)\le\sum_{i=0}^a\ind_{\{U_i=v\}}
          +\sum_{i=0}^{n-a}\ind_{\{\widehat U_i=u-v\}}.
\]
If the sum exceeds $x$, one of the two counts exceeds $x/2$. Both
segments have length at most $\lceil n/2\rceil$, so
\eqref{eq:initial-density} and \eqref{eq:block-density} bound their
$\Pi_n^u$-masses by the corresponding unconditioned tails up to a
constant. This proves the first inequality of \eqref{eq:local-tails}.
It also gives $\|\ell_n(v)+1\|_{p,\mu}\le CpL_n$ by integration, since
$\mu(\Omega)\le C$.

All sites visited before $n$ lie in $[-n,n]^2$. Taking a union bound
gives the second inequality of \eqref{eq:local-tails}. To get the
moment estimate, choose $x_0=C_0L_n^2$ large enough that
$(n+1)^2e^{-cx_0/L_n}\le1$. Then
$\mu(\ell_n^*>x_0+y)\le Ce^{-cy/L_n}$ for $y\ge0$.
The $L^p(\mu)$ norm of $(\ell_n^*-x_0)_+$ is therefore at most
$CpL_n$, giving the claimed bound for $\ell_n^*$.

It remains to count visits to the random site $U_j$ at a fixed index
$j$. Under $\PP$, visits before $j$ correspond to zeros of the path
with increments $Y_j,Y_{j-1},\ldots,Y_1$. Visits from $j$ onwards
correspond to zeros of the path with increments $Y_{j+1},\ldots,Y_n$.
Each has the law of an ordinary planar walk of length at most $n$, so
the fixed-site estimate and the triangle inequality give the desired
moment. For a bridge, let $V$ be the path started at zero with increments
$Y_{j+1},\ldots,Y_n,Y_1,\ldots,Y_j$. If $j\le i<n$ and $U_i=U_j$,
then $V_{i-j}=0$. If $0\le i<j$ and $U_i=U_j$, then
$V_{n-j+i}=u-U_j+U_i=u$. Both index ranges lie in $[0,n)$.
Thus $\ell_n(U_j)$ is bounded by the total local times of $V$ at $0$
and $u$ before $n$. This cyclic permutation preserves the joint law
of the planar increments under the bridge.
Applying the fixed-site bounds at $0$ and $u$ proves
\eqref{eq:local-index}, also when $j=0$ or $j=n$.
\end{proof}

\begin{lemma}\label{lem-range}
Put $b_s=\PP(U_i\ne0,\ 1\le i\le s)$, with $b_0=1$. Then
$b_s\le C/\log(s+2)$. For $r\ge1$ and $p\ge1$,
\begin{equation}\label{eq:range-unconditioned}
 \big\|\,|\{U_0,\ldots,U_r\}|\,\big\|_{p,\PP}
     \le C_p\,r/\log(r+2).
\end{equation}
In \eqref{eq:range-unconditioned} the constant is independent of $p$
when $1\le p\le\sqrt r$.
If $r=n-m$ and $n/2\le m<n$, then, uniformly in $u$,
\begin{equation}\label{eq:range-bridge}
 \big\|\,|\{U_m,\ldots,U_n\}|\,\big\|_{p,\Pi_n^u}
     \le C_p\,r/\log(r+2).
\end{equation}
\end{lemma}

\begin{proof}
The probability that the last visit to $0$ up to time $s$ occurs at
$i$ is $p_i(0)b_{s-i}$. Since $b_s$ is decreasing, summation gives
\[
 1=\sum_{i=0}^s p_i(0)b_{s-i}\ge b_s\sum_{i=0}^s p_i(0).
\]
The lower bound following \eqref{eq:planar-kernel-uniform} proves the
claim about $b_s$. Write $R=|\{U_0,\ldots,U_r\}|$. The site $U_i$ is
new precisely when the reversed partial sums $U_i-U_{i-k}$ are
nonzero for all $1\le k\le i$. By independence of the increments this
event has probability $b_i$. Thus
$\EE R=\sum_{i=0}^r b_i\le Cr/\log(r+2)$, by splitting the sum at
$\lfloor\sqrt r\rfloor$.

Given the path to any index $i$, the number of new sites still to be
visited is at most $|\{U_i,\ldots,U_r\}|$. By translation and independence
the conditional expectation of this last quantity is the expected range
of a walk of length $r-i$, hence at most $\EE R$. Apply the counting
argument of Lemma \ref{lem:bridge-local} to the first-visit indicators,
with $A=\EE R$. We obtain
\[
 \EE\exp\!\left(c R\log(r+2)/r\right)\le 2,
\]
which implies \eqref{eq:range-unconditioned} for each fixed $p$.

To make the constant uniform when $1\le p\le\sqrt r$, divide the $r$
increments into $a=\lceil p\rceil$ consecutive blocks, each of length
at most $b=\lceil r/a\rceil$. The range size of each block depends
only on its own increments, because translating the starting point
preserves the number of sites. These sizes are independent and their
sum dominates $R$. Each block has range stochastically bounded by that
of a walk of length $b$. Multiplying the preceding exponential moment
bounds, now at length $b$, gives
\[
 \EE\exp\!\left(cR\log(b+2)/b\right)\le2^a,\qquad
 \|R\|_{p,\PP}\le C\,2^{a/p}\frac{pb}{\log(b+2)}.
\]
The second inequality uses $y^p\le(p/e)^p e^y$. Here $a\le2p$,
$pb\le r+p\le2r$ and $b\ge r/(2p)\ge\sqrt r/2$, so
$\log(b+2)\ge c\log(r+2)$. This proves the uniform bound.
On $\{U_n=u\}$ the terminal range has exactly the size of the range of
$\widehat U_0,\ldots,\widehat U_r$. Since $r\le n/2$,
\eqref{eq:block-density} compares its $p$th moment with the
unconditioned one up to a fixed factor. Taking a $p$th root proves
\eqref{eq:range-bridge}.
\end{proof}

Conditioning on $U_n=u$ can change the law of a first departure, especially
near the last index. The walk can enter or leave column $v$ only on a
planar first departure from that column or one of its four neighboring
columns. We estimate the accumulated change in direction probabilities
at first departures from these five columns.

\begin{lemma}\label{lem:bridge-bias}
For $0\le j<n$ and $e$ in the set of five increments in
\eqref{eq:independent-drivers}, set
\[
 q_j(e)=Q(Y_{j+1}=e\mid\mathcal D_j),\qquad
 d_j=\max_e|q_j(e)-\nu_\theta(e)|.
\]
For every integer $1\le a\le n/2$ and $x\ge0$,
\begin{equation}\label{eq:direction-max}
 \mu\!\left(\max_{j\le n-a}d_j>x\right)\le C e^{-2ax^2}.
\end{equation}
Write $N[v]=\{v\}\cup\{w:w\sim v\}$ for a site and its four planar
neighbors. For each $p\ge1$,
\begin{equation}\label{eq:column-bias}
 J_n(v):=\sum_{j<n}\ind_{\{U_j\in N[v]\}}d_j,\qquad
 \sup_v\|J_n(v)\|_{p,\mu}\le Cp^{3/2} .
\end{equation}
Here and below a sum over $j<n$ starts at $j=0$.
\end{lemma}

\begin{proof}
Under $\PP$ the next increment has law $\nu_\theta$ given
$\mathcal D_j$, so every $d_j$ is zero. Take $Q=Q_n^u$ and
$\mu=\Pi_n^u$. Given $\mathcal D_j$, the unused increments
$Y_{j+1},\ldots,Y_n$ have their product law conditioned on their sum
being $u-U_j$. Indeed, $\mathcal D_j$ reveals $Y_1,\ldots,Y_j$, and
the additional history through $\tau_j$ is determined by those
increments and the independent walk $\widetilde H$. The conditional law
of the remaining increments is therefore invariant under permutations.
In particular their conditional probabilities of being equal to a
fixed $e$ are all the same. It follows that
\[
 q_j(e)=\EE_Q[\ind_{\{Y_n=e\}}\mid\mathcal D_j],\qquad 0\le j<n.
\]
Thus $q_j(e)$ is a martingale as $j$ varies. If $B_a$ is the frequency
of $e$ among the last $a$ increments, the same symmetry gives
\[
 B_a=\frac1a\sum_{i=n-a+1}^n\ind_{\{Y_i=e\}},\qquad
 q_{n-a}(e)=\EE_Q[B_a\mid\mathcal D_{n-a}].
\]

Write $\nu=\nu_\theta(e)$. Conditional Jensen's inequality under $Q$, followed by the tower property
under $Q$ and the identity $\mu=w_{n,u}Q$, gives the first inequality below. Applying
\eqref{eq:block-density} to the last $a$ increments gives the second.
For every $t\in\RR$,
\[
 \int e^{t(q_{n-a}(e)-\nu)}\,d\mu
 \le \int e^{t(B_a-\nu)}\,d\mu
 \le C\,\EE\exp\!\left(\frac ta\sum_{i=1}^a
                   (\ind_{\{Y_i=e\}}-\nu)\right)
 \le C e^{t^2/(8a)}.
\]
Under $\PP$ the indicators in the sum are independent Bernoulli
variables of mean $\nu$. Hoeffding's lemma bounds each centered
exponential moment by $\exp(t^2/(8a^2))$, giving the last inequality.
For $t>0$, the process $\exp(t(q_j(e)-\nu))$ is a nonnegative
submartingale under $Q$. Doob's maximal inequality and $\mu=w_{n,u}Q$ yield
\[
 \mu\!\left(\max_{0\le j\le n-a}(q_j(e)-\nu)>x\right)
       \le C\exp(-tx+t^2/(8a)).
\]
Take $t=4ax$ for $x>0$. Applying the same argument to $\nu-q_j(e)$
and summing over the five choices of $e$ yields \eqref{eq:direction-max}.
At $x=0$ the bound follows from $\mu(\Omega)\le C$. Integration of
this Gaussian tail gives, for every $p\ge1$,
\[
 \Big\|\max_{0\le j\le n-a}d_j\Big\|_{p,\mu}\le C\sqrt{p/a}.
\]

We now prove \eqref{eq:column-bias}. For $n=1$ the single possible summand in $J_n(v)$ is at most one and $\mu(\Omega)\le C$. Now suppose that $n\ge2$, and put $m=\lfloor n/2\rfloor$. At indices $j$ with $n-j>m$,
the maximum deviation has $L^{2p}(\mu)$ norm at most $C\sqrt{p/m}$.
The number of visits is at most $\sum_{w\in N[v]}\ell_n(w)$, whose
$L^{2p}(\mu)$ norm is at most $CpL_n$ by Lemma \ref{lem:bridge-local}.
By H\"older's inequality, 
\[
 \left\|\sum_{j< n-m}\ind_{\{U_j\in N[v]\}}d_j\right\|_{p,\mu}
 \le \Big\|\sum_{w\in N[v]}\ell_n(w)\Big\|_{2p,\mu}\,
       \Big\|\max_{j< n-m}d_j\Big\|_{2p,\mu}
 \le Cp^{3/2} L_n/\sqrt m,
\]
where the last term is bounded by $Cp^{3/2}$ uniformly in $n$.

For the remaining indices and each $h\ge0$, put
\[
 \mathcal J_h=\{j:1\le n-j\le m,\ 2^h\le n-j<2^{h+1}\},\qquad
 N_h(v)=\sum_{j\in \mathcal J_h}\ind_{\{U_j\in N[v]\}} .
\]
If $\mathcal J_h$ is nonempty, then $2^h\le m$ and $j\le n-2^h$ throughout
$\mathcal J_h$. Taking $a=2^h$ in the maximal estimate gives
$\|\max_{j\in \mathcal J_h}d_j\|_{2p,\mu}\le C\sqrt p\,2^{-h/2}$.
To bound $N_h(v)$, read the path backwards $\widehat U_i=U_n-U_{n-i}$ for $i\leq b=\min(2^{h+1},m)$. Every visit counted by $N_h(v)$ is a visit
by $\widehat U$ through index $b$ to one of the five sites $u-w$,
$w\in N[v]$. By \eqref{eq:block-density} and the unconditioned
fixed-site estimate, including the endpoint, we obtain
$\|N_h(v)\|_{2p,\mu}\le Cp\log(b+2)\le Cp(h+1)$. Consequently
\[
 \left\|\sum_{j\in \mathcal J_h}\ind_{\{U_j\in N[v]\}}d_j\right\|_{p,\mu}
 \le \|N_h(v)\|_{2p,\mu}\,
       \Big\|\max_{j\in \mathcal J_h}d_j\Big\|_{2p,\mu}
 \le Cp^{3/2}(h+1)2^{-h/2}.
\]
The series over $h$ converges, so the triangle inequality proves
\eqref{eq:column-bias}. 
\end{proof}

\section{Vertical occupation estimates}\label{sec:column}

At physical time $t$, write $P_t\in\ZZ^2$ for the planar position of
the walk and $Z_t\in\ZZ$ for its height, so that $S_t=(P_t,Z_t)$.
At discovery time $\tau_j$ these coordinates are $(U_j,H_j)$. Put
$\epsilon_t=Z_{t+1}-Z_t$, which equals $\pm1$ on a vertical step and
zero on a planar step. Excursion $j$ makes all its $D_j$ vertical
steps in column $U_j$. Thus the vertical time spent in column $v$
before $\tau_n$ has the following two expressions.
\[
 T_n(v):=\sum_{j=1}^n D_j\ind_{\{U_j=v\}}
        =\sum_{0\le t<\tau_n}\ind_{\{P_t=v\}}\epsilon_t^2,
 \qquad T_n=\sum_{v\in\ZZ^2} T_n(v).
\]
We first estimate $T_n(v)$ for a fixed column $v$, then sum over the
columns reached during a window of discovery indices.

\subsection{An occupation identity}

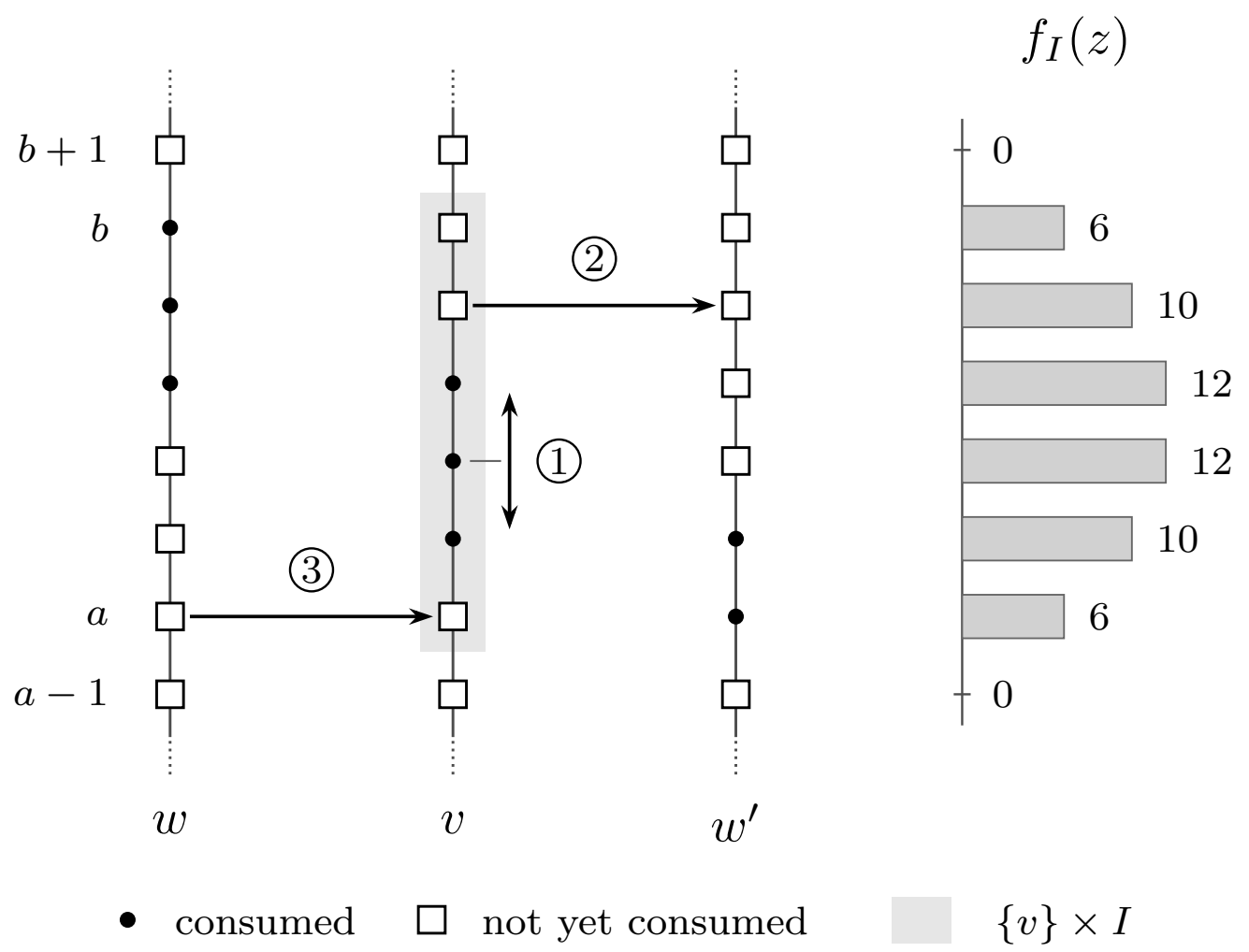
\begin{figure}[tbp]
\centering
\begin{tikzpicture}[
    x=1cm,y=0.55cm,
    >=Stealth,
    font=\small,
    consumed/.style={circle,fill=black,inner sep=0pt,minimum size=3.2pt},
    unvisited/.style={rectangle,draw=black,fill=white,line width=0.55pt,
                     inner sep=0pt,minimum size=5.4pt},
    move/.style={->,line width=0.7pt,
                 >={Stealth[length=1.7mm,width=1.15mm]}},
    tag/.style={circle,draw=black,line width=0.45pt,inner sep=0.7pt,
                font=\scriptsize},
    lab/.style={font=\scriptsize}
]
\fill[black!10] (3.27,0.55) rectangle (3.73,6.45);
\foreach \x in {1.5,3.5,5.5}{
  \draw[black!70,line width=0.55pt] (\x,-0.55) -- (\x,7.55);
  \draw[densely dotted,black!70,line width=0.55pt] (\x,7.6) -- (\x,8.05);
  \draw[densely dotted,black!70,line width=0.55pt] (\x,-0.6) -- (\x,-1.05);
}
\foreach \h in {4,5,6}{\node[consumed] at (1.5,\h) {};}
\foreach \h in {0,1,2,3,7}{\node[unvisited] at (1.5,\h) {};}
\foreach \h in {2,3,4}{\node[consumed] at (3.5,\h) {};}
\foreach \h in {0,1,5,6,7}{\node[unvisited] at (3.5,\h) {};}
\foreach \h in {1,2}{\node[consumed] at (5.5,\h) {};}
\foreach \h in {0,3,4,5,6,7}{\node[unvisited] at (5.5,\h) {};}
\foreach \h/\lab in {0/$a-1$,1/$a$,6/$b$,7/$b+1$}{
  \node[anchor=east,lab] at (1.2,\h) {\lab};
}
\node at (1.5,-1.65) {$w$};
\node at (3.5,-1.65) {$v$};
\node at (5.5,-1.65) {$w'$};

\draw[<->,line width=0.7pt,>={Stealth[length=1.7mm,width=1.15mm]}]
  (3.9,2.12) -- (3.9,3.88);
\draw[black!60,line width=0.4pt] (3.62,3) -- (3.84,3);
\node[tag] at (4.25,3.0) {1};
\draw[move] (3.64,5) -- (5.36,5);
\node[tag] at (4.5,5.6) {2};
\draw[move] (1.64,1) -- (3.36,1);
\node[tag] at (2.5,1.6) {3};

\def\xf{7.1}
\draw[black!70,line width=0.45pt] (\xf,-0.4) -- (\xf,7.4);
\foreach \h/\f in {1/6,2/10,3/12,4/12,5/10,6/6}{
  \fill[black!18] (\xf,\h-0.28) rectangle (\xf+0.12*\f,\h+0.28);
  \draw[black!60,line width=0.35pt] (\xf,\h-0.28) rectangle (\xf+0.12*\f,\h+0.28);
  \node[lab,anchor=west] at (\xf+0.12*\f+0.04,\h) {$\f$};
}
\foreach \h in {0,7}{
  \draw[black!70,line width=0.45pt] (\xf-0.06,\h) -- (\xf+0.06,\h);
  \node[lab,anchor=west] at (\xf+0.08,\h) {$0$};
}
\node at (\xf+0.8,8.4) {$f_I(z)$};

\node[consumed] at (1.2,-2.9) {};
\node[anchor=west,font=\scriptsize,text height=1.5ex,text depth=.25ex] at (1.4,-2.9) {consumed};
\node[unvisited] at (3.35,-2.9) {};
\node[anchor=west,font=\scriptsize,text height=1.5ex,text depth=.25ex] at (3.57,-2.9) {not yet consumed};
\fill[black!10] (6.6,-3.2) rectangle (7.02,-2.6);
\node[anchor=west,font=\scriptsize,text height=1.5ex,text depth=.25ex] at (7.2,-2.9) {$\{v\}\times I$};
\end{tikzpicture}
\caption{An interval of six heights in column $v$, two of its four
neighboring columns, and the exit-time function $f_I$. A vertical departure from $\{v\}\times I$ (1) increases $T_I$ by one. First departures from
the central column (2) and its neighbors (3) change the correction
$B_t(v,I)$ in \eqref{eq:boundary}.}
\label{fig:column-identity}
\end{figure}

Fix $v\in\ZZ^2$ and an integer interval $I=[a,b]\cap\ZZ$. Put
\begin{equation}
 \begin{split}
 f_I(z)&=(z-a+1)(b+1-z)\ind_{\{z\in I\}},\\
 h_I(z)&=(a+b-2z)\ind_{\{z\in I\}},\qquad
 g_I(x,z)=\ind_{\{x=v\}}f_I(z).
 \end{split}
\end{equation}
For $z\in I$, $f_I(z)$ is the mean exit time from $I$ for a
one-dimensional simple random walk started at $z$. A vertical step from
$\{v\}\times I$ therefore decreases $g_I$ by one on average. To use
this decrease to count vertical steps, we must also account for planar
steps into and out of the column. These occur only on first departures.
To compensate for their expected contribution under $\PP$, set
\begin{equation}\label{eq:boundary}
 B_t(v,I)=\theta\sum_{z\in I} f_I(z)
 \left(\ind_{\mathcal R_t}(v,z)-
              \frac14\sum_{w\sim v}\ind_{\mathcal R_t}(w,z)\right),
\end{equation}
where $w\sim v$ means that $w$ is a planar neighbor of $v$.
Figure \ref{fig:column-identity} shows the sites that contribute to
this correction.

After the first departure from $W_j$, the walk stays in the same column
until it reaches $W_{j+1}$. Thus $U_{j+1}=P_{\tau_j+1}$ is already known
at time $\tau_j+1$. Given the history at $\tau_j$, the probability under
$\PP$ that $U_{j+1}=v$ is $\nu_\theta(v-U_j)$. Write
\begin{equation}\label{eq:first-coefficient}
 C_j(v)=\ind_{\{U_{j+1}=v\}}-\nu_\theta(v-U_j)
\end{equation}
for the centered indicator.

\begin{lemma}\label{lem:column-identity}
Let $\zeta_I=\inf\{t\ge0:P_t=v,\ Z_t\in\{a-1,b+1\}\}$ and let
\[
 T_I(t)=\sum_{0\le r<t}\ind_{\{P_r=v,\,Z_r\in I\}}\epsilon_r^2
\]
count vertical departures from $\{v\}\times I$ before time $t$.
For $0\le t\le\zeta_I$ we have the pathwise identity
\begin{equation}\label{eq:column-increments}
 T_I(t)+g_I(S_t)-g_I(S_0)+B_t(v,I)
 =\sum_{0\le\tau_j<t} f_I(H_j)C_j(v)+\sum_{0\le r<t}\ind_{\{P_r=v\}}h_I(Z_r)\epsilon_r .
\end{equation}
Under $\PP$ the right-hand side, stopped at $\zeta_I$, is a martingale
starting from zero.
\end{lemma}

\begin{proof}
On a vertical step from $(v,z)$ with $z\in I$, the change in $g_I$
is $f_I(z+\epsilon_r)-f_I(z)=-1+\epsilon_r h_I(z)$. This also holds
on a step out of $I$, since the quadratic defining $f_I$ vanishes at
$a-1$ and $b+1$. The increment of $T_I$ cancels the $-1$.
A vertical step entering $\{v\}\times I$ from outside would have to
depart from one of the two sites defining $\zeta_I$, so no such step
occurs before the stop.

At a first departure from $W_j$, the change of column contributes
\[
 f_I(H_j)\bigl(\ind_{\{U_{j+1}=v\}}-\ind_{\{U_j=v\}}\bigr)
\]
to the increment of $g_I$. The departure also adds $W_j$ to the
consumed set, even when the step is vertical and the change of column
contributes zero. By \eqref{eq:boundary}, the increment of $B_t(v,I)$ is
\[
 f_I(H_j)\theta\left(\ind_{\{U_j=v\}}-\frac14\ind_{\{U_j\sim v\}}\right)
 =f_I(H_j)\bigl(\ind_{\{U_j=v\}}-\nu_\theta(v-U_j)\bigr).
\]
The two displayed terms sum to $f_I(H_j)C_j(v)$. Summing all increments
up to time $t\le\zeta_I$ proves \eqref{eq:column-increments}.

Under $\PP$, $H_j$ is known at time $\tau_j$ and $C_j(v)$ has
conditional mean zero given $\FF_{\tau_j}$. The increments of the
second sum also have conditional mean zero since each vertical sign is
symmetric given the history and the choice of a vertical step. This
proves the martingale assertion.
\end{proof}

\subsection{Occupation of a column}

Let $\mathcal I_m(v)$ be the collection of maximal integer intervals
of consumed heights in column $v$ at discovery time $\tau_m$. If
$I=[a,b]\cap\ZZ$ is one of these intervals, the two sites
$(v,a-1)$ and $(v,b+1)$ are outside $\mathcal R_{\tau_m}$. The walk
therefore did not visit either site strictly before $\tau_m$. It may
reach one of them at $\tau_m$ itself. Thus \eqref{eq:column-increments}
holds up to $\tau_m$ for this interval.
Moreover, $\mathbf1_{\mathcal R_{\tau_m}}(v,z)=1$ for every $z\in I$, and the current vertex $W_m$ is fresh. Hence
\[
 B_{\tau_m}(v,I)\ge0,\qquad g_I(S_{\tau_m})=0.
\]

The intervals in $\mathcal I_n(v)$ are determined by the path up to
$\tau_n$. If their endpoints were used directly as coefficients of
martingale increments, those coefficients could depend on future
steps. We will instead express the two sums in \eqref{eq:column-increments}
by Fourier inversion. For each fixed frequency, the same two random
sums serve every interval and can be estimated before the final
components are chosen. We note that \eqref{eq:column-work} below will suffice for transience, and \eqref{eq:column-uniform-orders} will be used to obtain the \((\log\log t)^3\) factor in the almost sure range bound of Proposition \ref{rangelowerbd}.

\begin{lemma}\label{lem:column-occupation}
For $p\ge1$ and either choice of $\mu$ in Section \ref{sec:bridge},
\begin{equation}\label{eq:column-work}
 \sup_{v\in\ZZ^2}\|T_n(v)\|_{p,\mu}\le C_p L_n^{5/2}.
\end{equation}
In particular, uniformly in the endpoint $y$ and column $v$,
\begin{equation}\label{eq:weighted-column}
 (n+1)\EE[T_n(v)^p\ind_{\{U_n=y\}}]\le C_p L_n^{5p/2}.
\end{equation}
More precisely, for $\mathcal L\ge\max\{1,L_n,p^4\}$,
\begin{equation}\label{eq:column-uniform-orders}
 \sup_v\|T_n(v)\|_{p,\mu}\le C p^3 \mathcal L^{5/2},
\end{equation}
with $C$ depending only on $\theta$.
\end{lemma}

\begin{proof}
Fix a column $v$ and put $\ell=\ell_n(v)$. We prove
\eqref{eq:column-uniform-orders} for $\mathcal L\ge\max\{1,L_n,p^4\}$.
Write $X=\|T_n(v)\|_{p,\mu}$. This is finite for each $n$, since
$T_n(v)\le\sum_{j=1}^n D_j$ and \eqref{eq:crude-exit} applies with
$k=j$ to each summand.

We first rewrite the two sums in \eqref{eq:column-increments} using
random sums that do not depend on the final interval decomposition. For $\xi\in[-\pi,\pi]$, put
\[
 F(\xi)=\sum_{j<n} e^{i\xi H_j} C_j(v)
 \qquad\text{and}\qquad
 M(\xi)=\sum_{r<\tau_n} e^{i\xi Z_r}\ind_{\{P_r=v\}}\epsilon_r .
\]
With the convention $\widehat f(\xi)=\sum_z f(z)e^{-i\xi z}$, Fourier
inversion gives, for every integer interval $I$,
\[
 \begin{aligned}
 \sum_{j<n} f_I(H_j) C_j(v)
 &=\frac1{2\pi}\int_{-\pi}^\pi \widehat f_I(\xi) F(\xi)\,d\xi,\\
 \sum_{r<\tau_n}\ind_{\{P_r=v\}}h_I(Z_r)\epsilon_r
 &=\frac1{2\pi}\int_{-\pi}^\pi \widehat h_I(\xi) M(\xi)\,d\xi .
 \end{aligned}
\]
We now sum the occupation identity over the final components
$I\in\mathcal I_n(v)$. As observed above, for each such $I$ the
identity holds up to $\tau_n$, with $B_{\tau_n}(v,I)\ge0$ and
$g_I(S_{\tau_n})=0$. Every vertical departure from column $v$
before $\tau_n$ starts in exactly one of these intervals, so
$\sum_I T_I(\tau_n)=T_n(v)$. At most one interval contains the
starting vertex in column $v$, and its contribution to the initial
term is at most $\ell^2$. The inversion formulas therefore give
\begin{equation}
    \label{eq:Tnvinver}
     T_n(v)\le\ell^2+\frac1{2\pi}\int_{-\pi}^{\pi}
 \left[\left(\sum_I|\widehat f_I(\xi)|\right)|F(\xi)|
       +\left(\sum_I|\widehat h_I(\xi)|\right)|M(\xi)|\right]\,d\xi .
\end{equation}

To estimate $F(\xi)$, let $m_j(v)=\EE_Q[C_j(v)\mid\mathcal D_j]$.
For $Q=\PP$ this mean is zero. To control it under the planar bridge,
recall from Lemma \ref{lem:bridge-bias} that $N[v]$
consists of $v$ and its four planar neighbors. For $U_j\in N[v]$ we
have $m_j(v)=q_j(v-U_j)-\nu_\theta(v-U_j)$, and for $U_j\notin N[v]$
both terms in \eqref{eq:first-coefficient} vanish. Thus
\begin{equation}\label{eq:coefficient-drift}
 |m_j(v)|\le\ind_{\{U_j\in N[v]\}}d_j .
\end{equation}
Fix $\xi$ and let $A_k$ be the real part of the centered sum up to
discovery $k$:
\[
 A_k=\Re\sum_{j=0}^{k-1}e^{i\xi H_j}\bigl(C_j(v)-m_j(v)\bigr),
 \qquad 0\le k\le n.
\]
Since $H_j$ is $\mathcal D_j$-measurable and
\[
 C_j(v)-m_j(v)
 =\ind_{\{U_{j+1}=v\}}-Q(U_{j+1}=v\mid\mathcal D_j),
\]
$(A_k)$ is a martingale under $Q$. Its $j$th increment satisfies
$|A_{j+1}-A_j|\le b_j$, where $b_j=\ind_{\{U_j\in N[v]\}}$
is known at discovery $j$. The sum of the squared increment bounds
therefore counts the visits of $U$ to $N[v]$ before index $n$:
\[
 \sum_{j=0}^{n-1}b_j^2
 =\sum_{j=0}^{n-1}\ind_{\{U_j\in N[v]\}}
 =\sum_{w\in N[v]}\ell_n(w)=:N_v .
\]
By conditional Hoeffding's lemma,
\[
 \EE_Q\!\left[e^{s(A_{j+1}-A_j)}\mid\mathcal D_j\right]
 \le e^{s^2b_j^2/2}.
\]
Iterating gives $\EE_Q[e^{sA_n-s^2N_v/2}]\le1$. Applying this bound
to $A_n$ and $-A_n$ with $s=x/y$ gives, for $x,y>0$,
\[
 Q(|A_n|>x,\,N_v\le y)\le 2e^{-x^2/(2y)}.
\]
Since $\mu$ is a bounded multiple of $Q$, it follows that
\[
 \mu(|A_n|>x)
 \le C e^{-x^2/(2y)}+\mu(N_v>y)
 \le C e^{-x^2/(2y)}+C e^{-cy/L_n}.
\]
The last inequality follows from \eqref{eq:local-tails}, since
$N_v>y$ implies $\ell_n(w)>y/5$ for some $w\in N[v]$.
 Taking $y=x\sqrt{L_n}$ in the preceding estimate gives
\[
 \mu(|A_n|>x)\le Ce^{-cx/\sqrt{L_n}} .
\]
Consequently,
\[
\begin{aligned}
 \int |A_n|^{2p}\,d\mu
 &=2p\int_0^\infty x^{2p-1}\mu(|A_n|>x)\,dx\\
 &\le Cp\int_0^\infty x^{2p-1}e^{-cx/\sqrt{L_n}}\,dx\\
 &=Cp\,\Gamma(2p)(c^{-1}\sqrt{L_n})^{2p}
 \le (Cp\sqrt{L_n})^{2p},
\end{aligned}
\]
where the last step uses $\Gamma(2p)^{1/(2p)}\le Cp$.
Taking the $2p$th root gives
$\|A_n\|_{2p,\mu}\le Cp\sqrt{L_n}$.
The imaginary part satisfies the same bound. The absolute value of
the drift $\sum_{j<n}e^{i\xi H_j}m_j(v)$ is at most $J_n(v)$ by
\eqref{eq:coefficient-drift}. Its $L^{2p}(\mu)$ norm is at most
$Cp^{3/2}$ by \eqref{eq:column-bias}. Hence
$\|F(\xi)\|_{2p,\mu}\le Cp\sqrt{\mathcal L}$.

To estimate $M(\xi)$, we may reveal the entire planar path $U$.
In the filtration $\GG_r=\FF_r\vee\sigma(U)$, we then know whether
the next step is vertical, while its sign on a vertical step is
still symmetric by Lemma~\ref{lem:height-martingale}(a).
Thus the partial sums defining $M(\xi)$, stopped at $\tau_n$,
form a martingale under $Q$. Their squared absolute increments sum to
\[
 \sum_{r<\tau_n}
 \bigl|e^{i\xi Z_r}\ind_{\{P_r=v\}}\epsilon_r\bigr|^2=T_n(v).
\]
The $L^{2p}(\mu)$ norm of its square root is therefore
$\|T_n(v)^{1/2}\|_{2p,\mu}=\sqrt X$.

We use the following martingale moment inequality, with explicit
dependence on the moment order. Let $\mathcal M$ be a real martingale
with $\mathcal M_0=0$ and bounded increments under $Q$. Write
$\mathcal M_b^*=\max_{j\le b}|\mathcal M_j|$ and
$[\mathcal M]_b=\sum_{j=1}^b(\mathcal M_j-\mathcal M_{j-1})^2$.
For $q\ge2$, we apply Taylor's formula to $|\mathcal M_j|^q$ about
$\mathcal M_{j-1}$. The linear terms have zero conditional mean, and
the second derivative at the intermediate point is at most
$q(q-1)(\mathcal M_b^*)^{q-2}$. Summing over $j$ and then applying
H\"older's inequality with exponents $q/(q-2)$ and $q/2$ gives
\[
 \EE_Q|\mathcal M_b|^q
 \le \frac12q(q-1)\,
 \EE_Q\!\left[(\mathcal M_b^*)^{q-2}[\mathcal M]_b\right] \le \frac12q(q-1)\,
 \|\mathcal M_b^*\|_{q,Q}^{q-2}
 \|[\mathcal M]_b^{1/2}\|_{q,Q}^2 .
\]
Doob's inequality and $(q/(q-1))^{q-2}\le e$ now give
$\|\mathcal M_b\|_{q,Q}
 \le 2q\,\|[\mathcal M]_b^{1/2}\|_{q,Q}$
(see also \cite[Theorem 1.1]{BDG1972}). The inequality holds with
$\mu$ in both norms by homogeneity. Apply it with $q=2p$ to the real
and imaginary parts of the vertical sum stopped at $\tau_n\wedge b$. For each of these stopped martingales, the sum of the squared increments is at most
$T_n(v)$. Fatou's lemma as
$b\to\infty$, together with the bound from the preceding step, yields
\begin{equation}\label{eq:frequency-moments}
 \sup_\xi\|F(\xi)\|_{2p,\mu}\le Cp\sqrt{\mathcal L}
 \qquad\text{and}\qquad
 \sup_\xi\|M(\xi)\|_{2p,\mu}\le Cp\sqrt X .
\end{equation}

It remains to bound
$\sum_{I\in\mathcal I_n(v)}|\widehat f_I(\xi)|$ and
$\sum_{I\in\mathcal I_n(v)}|\widehat h_I(\xi)|$
in terms of $\ell$, the number of consumed vertices in column $v$.
Fix $I=[a,b]\cap\ZZ$ and write $k=|I|=b-a+1$. For every $\xi$,
the triangle inequality gives the absolute-sum bounds
\begin{equation}\label{eq:fourier-absolute-sums}
 |\widehat f_I(\xi)|\le\|f_I\|_1\le k^3
 \qquad\text{and}\qquad
 |\widehat h_I(\xi)|\le\|h_I\|_1\le k^2 .
\end{equation}
Define the forward difference
of $h_I$ and the centered second difference of $f_I$ by
\[
 \Delta h_I(z)=h_I(z+1)-h_I(z),\qquad
 \Delta^2f_I(z)=f_I(z+1)+f_I(z-1)-2f_I(z).
\]
Direct computation gives
\begin{equation}\label{eq:difference-ell-one}
 \|\Delta^2f_I\|_1=4k
 \qquad\text{and}\qquad
 \|\Delta h_I\|_1=4(k-1).
\end{equation}
Under our convention
$\widehat\phi(\xi)=\sum_z\phi(z)e^{-i\xi z}$, these differences
satisfy
\[
 \widehat{\Delta h_I}(\xi)
 =(e^{i\xi}-1)\widehat h_I(\xi),
 \qquad
 \widehat{\Delta^2f_I}(\xi)
 =-4\sin^2(\xi/2)\,\widehat f_I(\xi).
\]
Since the absolute Fourier transform of a finitely supported function
is at most its $\ell^1$ norm, \eqref{eq:difference-ell-one} yields
\[
 4\sin^2(\xi/2)\,|\widehat f_I(\xi)|
 \le4k,
 \qquad
 |1-e^{i\xi}|\,|\widehat h_I(\xi)|
 \le4(k-1).
\]
For $0<|\xi|\le\pi$, we have
\[
 4\sin^2(\xi/2)\ge\frac{4\xi^2}{\pi^2}
 \qquad\text{and}\qquad
 |1-e^{i\xi}|=2|\sin(\xi/2)|\ge\frac{2|\xi|}{\pi}.
\]
Therefore
\[
 |\widehat f_I(\xi)|\le\frac{Ck}{|\xi|^2}
 \qquad\text{and}\qquad
 |\widehat h_I(\xi)|\le\frac{Ck}{|\xi|}.
\]
Combining these decay bounds with
\eqref{eq:fourier-absolute-sums} gives
\begin{equation}\label{eq:fourier-single-interval}
 |\widehat f_I(\xi)|
 \le C\min\{k^3,k/|\xi|^2\}
 \qquad\text{and}\qquad
 |\widehat h_I(\xi)|
 \le C\min\{k^2,k/|\xi|\},
\end{equation}
with the second member of each minimum understood to be infinite at $\xi=0$. Since $\sum_I|I|=\ell$, we have
$\sum_I|I|^3\le\ell^3$ and $\sum_I|I|^2\le\ell^2$. Summing
\eqref{eq:fourier-single-interval} over
$I\in\mathcal I_n(v)$ therefore gives
\[
 \sum_I|\widehat f_I(\xi)|
 \le C\min\{\ell^3,\ell/|\xi|^2\},
 \qquad
 \sum_I|\widehat h_I(\xi)|
 \le C\min\{\ell^2,\ell/|\xi|\}.
\]
Substitution into \eqref{eq:Tnvinver} yields
\begin{equation}\label{eq:fourier-pathwise}
 T_n(v)\le\ell^2+C\int_{-\pi}^{\pi}
 \left[\min\{\ell^3,\ell/|\xi|^2\}|F(\xi)|
       +\min\{\ell^2,\ell/|\xi|\}|M(\xi)|\right]\,d\xi .
\end{equation}

We now take $L^p(\mu)$ norms in \eqref{eq:fourier-pathwise} to obtain
an inequality for $X$. Put $K=p\mathcal L$. For $d=1,2,3$,
Lemma~\ref{lem:bridge-local} gives
\[
 \|\ell^d\|_{2p,\mu}=\|\ell\|_{2dp,\mu}^d\le (CK)^d.
\]
By Minkowski's integral inequality, it suffices to bound the
$L^p(\mu)$ norm of each integrand. For the term involving $F(\xi)$,
H\"older's inequality gives, for $0<|\xi|\le\pi$,
\[
 \begin{aligned}
 \bigl\|\min\{\ell^3,\ell/|\xi|^2\}\,F(\xi)\bigr\|_{p,\mu}
 &\le \bigl\|\min\{\ell^3,\ell/|\xi|^2\}\bigr\|_{2p,\mu}
       \,\|F(\xi)\|_{2p,\mu} \\
 &\le C\min\{K^3,K/|\xi|^2\}\,\|F(\xi)\|_{2p,\mu}.
 \end{aligned}
\]
The same argument bounds the integrand involving $M(\xi)$ in
$L^p(\mu)$ by $C\min\{K^2,K/|\xi|\}\,\|M(\xi)\|_{2p,\mu}$.
To integrate the two bounds, we split at $1/K$, obtaining
\[
 \int_0^\pi \min\{K^3,K/\xi^2\}\,d\xi\le CK^2
 \qquad\text{and}\qquad
 \int_0^\pi \min\{K^2,K/\xi\}\,d\xi\le CK\log(2+K).
\]
Using the uniform bounds on $\|F(\xi)\|_{2p,\mu}$ and
$\|M(\xi)\|_{2p,\mu}$ from \eqref{eq:frequency-moments} and
substituting $K=p\mathcal L$, we obtain
\[
 X\le Cp^3\mathcal L^{5/2}
      +Cp^2\mathcal L\log(2+p\mathcal L)\sqrt X .
\]
As $x_1x_2\leq (x_1^2+x_2^2)/2$, the last term is upper bounded by
$X/2+Cp^4\mathcal L^2\log^2(2+p\mathcal L)$. Since
$p\le\mathcal L^{1/4}$, we have
$p\log^2(2+p\mathcal L)\le C\sqrt{\mathcal L}$, which proves
\eqref{eq:column-uniform-orders}. For fixed $p$, take
$\mathcal L=\max\{1,L_n,p^4\}$ to obtain \eqref{eq:column-work}.
Taking $\mu=\Pi_n^y$ and raising to the $p$th power gives
\eqref{eq:weighted-column}.
\end{proof}

\subsection{Vertical time in a terminal window}

\begin{proposition}\label{prop:window}
Let $p\ge1$, $0\le m<n$ and $r=n-m$. For $\mu=\PP$ or $\mu=\Pi_n^u$,
uniformly in $u$,
\begin{equation}\label{eq:window}
 \|T_n-T_m\|_{p,\mu}\le C_p\frac{rL_n^{5/2}}{\log(r+2)}.
\end{equation}
For a fixed admissible endpoint, division by $w_{n,u}^{1/p}$ gives
the same bound under $Q_n^u$, with a constant depending also on $u$.
For every fixed $\alpha>0$, this is $O(rL_n^{3/2})$ when $r\ge n^\alpha$,
with the constant also depending on $\alpha$. In particular,
\begin{equation}\label{eq:whole-clock}
 \|T_n\|_{p,\PP}\le C_p nL_n^{3/2}.
\end{equation}
For $1\le p\le L_n^{1/4}$ we also have
\begin{equation}\label{eq:clock-high-moments}
 \|T_n\|_{p,\PP}\le C p^3 n L_n^{3/2},
\end{equation}
with $C$ depending only on $\theta$.
\end{proposition}

\begin{proof}
Write $\Delta T=T_n-T_m$. Every vertical step counted in $\Delta T$
occurs in a column of $\mathcal V=\{U_m,\ldots,U_n\}$. Put
$F_v=\{v\in\mathcal V\}$ for $v\in\ZZ^2$.
H\"older's inequality on each path gives
\[
 (\Delta T)^p
 \le\left(\sum_vT_n(v)\ind_{F_v}\right)^p
 \le |\mathcal V|^{p-1}\sum_vT_n(v)^p\ind_{F_v}.
\]
Integrating and applying Cauchy--Schwarz with respect to the product
of $\mu$ and counting measure on $\ZZ^2$, we obtain
\begin{equation}\label{eq:selected-holder}
\begin{aligned}
 \int (\Delta T)^p\,d\mu
 &\le \sum_v \int
      \bigl(|\mathcal V|^{p-1}\ind_{F_v}\bigr)
      \bigl(T_n(v)^p\ind_{F_v}\bigr)\,d\mu \\
 &\le \left(\sum_v\int |\mathcal V|^{2p-2}\ind_{F_v}\,d\mu\right)^{1/2}
      \left(\sum_v\int T_n(v)^{2p}\ind_{F_v}\,d\mu\right)^{1/2} \\
 &= \left(\int |\mathcal V|^{2p-1}\,d\mu\right)^{1/2}
    \left(\sum_v\int T_n(v)^{2p}\ind_{F_v}\,d\mu\right)^{1/2},
\end{aligned}
\end{equation}
where the last equality uses
$\sum_v\ind_{F_v}=|\mathcal V|$. The first factor involves only the
planar range in the window. Lemma \ref{lem-range} gives
$$
 \bigl\|\,|\mathcal V|\,\bigr\|_{2p-1,\mu}
 \le C_p\frac r{\log(r+2)}
$$
under $\PP$ for every $m$ by stationarity of the increments, and under
$\Pi_n^u$ for $m\ge n/2$. We first work in these cases.

In the second factor, the columns whose occupation we count are
selected by the planar path during the window. To estimate their
contribution, we split $F_v$ according to the
last index $j\in[m,n]$ with $U_j=v$. If $j$ is the last such index,
then $T_n(v)=T_j(v)$, since no later excursion takes place in column
$v$. The last visit is not a stopping time. We therefore fix a
deterministic $j$ and condition on the history $\mathcal D_j$ under
$\PP$. On $\{U_j=v\}$ the remaining planar increments are independent
of this history. For $s\ge1$ put
\[
 a_s(x)=\PP(U_s=x,\ U_i\ne0\text{ for }1\le i\le s),\qquad
 a_0(x)=\ind_{\{x=0\}} .
\]
Given $\mathcal D_j$ on $\{U_j=v\}$, the probability that the path
avoids $v$ at indices $j+1,\ldots,n$ and ends at $u$ is $a_{n-j}(u-v)$. Summing over the possible endpoints
gives the no-return probability $b_{n-j}$ of Lemma \ref{lem-range}.
For a moment order $b\ge1$, write
$h_j(v)=\EE[T_j(v)^b\ind_{\{U_j=v\}}]$. Summing over the disjoint
possibilities for the last visit gives
\begin{align}
 \sum_v\int T_n(v)^b\ind_{F_v}\,d\Pi_n^u
   &=(n+1)\sum_{j=m}^n\sum_v h_j(v)a_{n-j}(u-v),
       \label{eq:last-visit}\\
 \sum_v\EE[T_n(v)^b\ind_{F_v}]
   &=\sum_{j=m}^n b_{n-j}\sum_v h_j(v). \label{eq:last-visit-free}
\end{align}
The term $j=n$ uses $a_0(u-v)=\ind_{\{u=v\}}$ and $b_0=1$.
The factor $n+1$ in \eqref{eq:last-visit} comes from the definition
of $\Pi_n^u$. For $m\ge n/2$, applying
\eqref{eq:weighted-column} at time $j$ with endpoint $v$ gives
$h_j(v)\le C_b L_n^{5b/2}/(j+1)$. Since $(n+1)/(j+1)\le2$ and
$\sum_v a_{n-j}(u-v)=b_{n-j}$,
\eqref{eq:last-visit} is at most $C_b L_n^{5b/2}\sum_{s=0}^r b_s$.

Under $\PP$, formula \eqref{eq:last-visit-free} requires a bound on
$\sum_v h_j(v)$. To retain the endpoint probability when summing over $v$,
we use Cauchy--Schwarz to write
\[
 h_j(v)\le
  \bigl(\EE[T_j(v)^{2b}\ind_{\{U_j=v\}}]\bigr)^{1/2}
  p_j(v)^{1/2}
  \le C_b L_n^{5b/2}\sqrt{p_j(v)/(j+1)} .
\]
The last inequality is \eqref{eq:weighted-column} with power $2b$.
To sum the square roots, put $g_j(v)=1+\|v\|_2^2/(j+1)$.
Independence and boundedness of the centered increments give
$\EE \|U_j\|_2^4\le C(j+1)^2$. The disk of radius $\sqrt{j+1}$ contains
$O(j+1)$ sites. On the annulus between radii $2^k\sqrt{j+1}$ and
$2^{k+1}\sqrt{j+1}$, $g_j^{-2}\le 2^{-4k}$ and there are
$O(2^{2k}(j+1))$ sites. Thus $\sum_v g_j(v)^{-2}\le C(j+1)$, and
\[
 \sum_v\sqrt{p_j(v)}
 \le \bigg(\sum_v g_j(v)^{-2}\bigg)^{1/2}
           \big(\EE g_j(U_j)^2\big)^{1/2}\le C\sqrt{j+1}.
\]
In particular, $\sum_v h_j(v)\le C_bL_n^{5b/2}$ uniformly in $j$.
In both last-visit formulas the remaining sum is $\sum_{s=0}^r b_s$.
This is the expected range of a planar walk of $r$ steps, since the
probability of discovering a new site at step $s$ is $b_s$ by reversal.
Lemma \ref{lem-range} therefore gives, under $\PP$ for all $m$ and
under $\Pi_n^u$ for $m\ge n/2$,
\begin{equation}\label{eq:selected-column}
 \sum_v\int T_n(v)^b\ind_{F_v}\,d\mu
     \le C_b L_n^{5b/2}\sum_{s=0}^r b_s
     \le C_b' L_n^{5b/2}\frac r{\log(r+2)} .
\end{equation}

Substituting \eqref{eq:selected-column} with $b=2p$ and the range
bound into \eqref{eq:selected-holder} gives
\[
 \int (\Delta T)^p\,d\mu\le
 C_p\left(\frac r{\log(r+2)}\right)^{p-1/2}
    L_n^{5p/2}\left(\frac r{\log(r+2)}\right)^{1/2}.
\]
Taking the $p$th root proves \eqref{eq:window} under $\PP$ for every
$m$ and under $\Pi_n^u$ for $m\ge n/2$. Taking $m=0$ under $\PP$
also gives \eqref{eq:whole-clock}.

To cover the remaining bridge windows, suppose $m<n/2$ under $\Pi_n^u$ and put
$a=\lceil n/2\rceil$. The clock $T_a$ is determined by the first $a$
planar increments and the independent vertical walk $\widetilde H$, so
\eqref{eq:initial-density} and \eqref{eq:whole-clock} give
$\|T_a\|_{p,\Pi_n^u}\le C_p nL_n^{3/2}$. Applying the already proved case $m\ge n/2$ with $m=a$ gives the same estimate for $T_n-T_a$ (which is zero if $a=n$).
Since $T_n-T_m\le T_a+(T_n-T_a)$, the resulting bound is
$C_p nL_n^{3/2}$. Here $r>n/2$ and $\log(r+2)\le L_n$, so this is
at most $2 C_p r L_n^{5/2}/\log(r+2)$, as claimed.

For \eqref{eq:clock-high-moments}, return to the two factors on the
right of \eqref{eq:selected-holder}, taking $m=0$, $\mu=\PP$ and
$b=2p$. We now track their dependence on $p$. The
Cauchy--Schwarz bound for $h_j(v)$ uses a column moment of order $4p$.
Because $p^4\le L_n$, the parameter $\mathcal L=4^4L_n$ is at least
$(4p)^4$ and $L_j$ for $j\le n$. Thus
\eqref{eq:column-uniform-orders}, applied under $\Pi_j^v$, gives
\[
 h_j(v)\le (C p^3 L_n^{5/2})^{2p}\sqrt{p_j(v)/(j+1)} .
\]
Summing over $v$ and using
$\sum_v\sqrt{p_j(v)}\le C\sqrt{j+1}$ give
\[
 \sum_v h_j(v)
 \le (Cp^3L_n^{5/2})^{2p}
 \sum_v\sqrt{\frac{p_j(v)}{j+1}}\le (Cp^3L_n^{5/2})^{2p},
\]
where the constant $C$ may change from line to line.
Substituting this into \eqref{eq:last-visit-free}, with $m=0$, gives
\[
\begin{aligned}
 \sum_v\EE\!\left[T_n(v)^{2p}\ind_{F_v}\right]
 &=\sum_{j=0}^n b_{n-j}\sum_v h_j(v) \le (Cp^3L_n^{5/2})^{2p}\sum_{s=0}^n b_s\le (Cp^3L_n^{5/2})^{2p}\frac{Cn}{L_n}.
\end{aligned}
\]
Thus the second factor in \eqref{eq:selected-holder} is at most  $(C p^3 L_n^{5/2})^p(Cn/L_n)^{1/2}$. The first factor is at most $(Cn/L_n)^{p-1/2}$ by the uniform range
bound \eqref{eq:range-unconditioned}, since $2p-1\le\sqrt n$ for large $n$.
This proves \eqref{eq:clock-high-moments}. For the finitely many
remaining $n$, apply \eqref{eq:whole-clock} at one fixed order above
all the allowed $p$ and use monotonicity of the norms under $\PP$.
\end{proof}

\section{A martingale defocusing estimate}
\label{sec:defocusing}

We adapt the space--time rectangle argument of Peres, Schapira and
Sousi \cite[Proposition 2.1]{MR3531697}. For a square-integrable real
martingale $(M_t,\FF_t)_{0\le t\le T}$, write
\[
 V_t:=\sum_{r=0}^{t-1}\EE[(M_{r+1}-M_r)^2\mid\FF_r]
\]
for its predictable quadratic variation, or bracket. The proposition
uses two estimates on different time intervals. On each block
$[s_i,s_{i+1}]$ we ask that the bracket increase by at least
$\lambda (T-s_i)$. Throughout the longer interval $[s_i,T]$ we require an upper
bound $B(T-s_i)$ for the conditional expectation of the remaining
bracket. The lower bound makes it unlikely that the martingale stays
in a narrow interval for the entire block. The upper bound controls
the chance of returning to that interval after the martingale has
moved sufficiently far away. A separate bound on the jumps controls
how far the martingale can overshoot an interval boundary. 

\begin{proposition}\label{lem:rectangle}
There exist absolute constants $\varepsilon,\eta,c>0$ with the following
property. Let $M$ be a square-integrable real martingale up to time $T$,
and let
\[
 0\le s_0<s_1<\cdots<s_N<T
\]
be deterministic times, where $N\ge1$. Set $k_i:=T-s_i$, and assume that $k_{i+1}\le k_i/2$ for $0\leq i<N$. Suppose that $0<\lambda\le B$
and
\begin{equation}\label{eq:abstract-jumps}
 |M_{t+1}-M_t|\le b\le\varepsilon\sqrt{\lambda k_N}
 \quad\text{a.s.},\qquad 0\le t<T.
\end{equation}
For every $0\le i<N$ and deterministic $s_i\le t\le T$, assume that
\begin{equation}\label{eq:rectangle-bracket}
 \EE[V_T-V_t\mid\FF_t]=\EE[(M_T-M_t)^2\mid\FF_t]\le Bk_i
 \quad\text{a.s.}
\end{equation}
Equivalently, the same bound holds at every stopping time
$\rho\in[s_i,T]$, by restriction to each event $\{\rho=t\}$. Put
\begin{equation}\label{eq:abstract-lower-events}
 A_i:=\{V_{s_{i+1}}-V_{s_i}\ge\lambda k_i\}.
\end{equation}
Then, for every $z\in\RR$,
\begin{equation}\label{eq:abstract-smallball}
 \PP\left(|M_T-z|\le\eta\sqrt{\lambda k_N},\,
          \bigcap_{i=0}^{N-1}A_i\right)
 \le\exp\left(-cN\sqrt{\lambda/B}\right).
\end{equation}
\end{proposition}

\begin{proof}
By replacing $M$ with $M-z$, it suffices to consider $z=0$.
Fix $\gamma:=1/10$, choose $\varepsilon\le\gamma/10$, and put
\[
 a_i:=\gamma\sqrt{\lambda k_i},\qquad
 \sigma_i:=\inf\{t\in\{s_i,\ldots,T\}:|M_t|\le a_i\},
\]
with $\inf\varnothing=\infty$. Thus $\sigma_i$ is the first entrance
into the rectangle $[s_i,T]\times[-a_i,a_i]$. Each rectangle contains
the next one, so $\sigma_i\le\sigma_{i+1}$, and a path that ends in
$[-a_N,a_N]$ must enter every rectangle (Figure \ref{fig:nested-rectangles}).
It suffices to show that upon each entrance there is a definite chance
of missing the next rectangle or of failing the lower bracket bound
on the current block. More precisely, we shall prove that on
$\{\sigma_i\le T\}$,
\begin{equation}\label{eq:rectangle-miss}
 \PP(\{\sigma_{i+1}>T\}\cup A_i^c\mid\FF_{\sigma_i\wedge T})
 \ge c_1\sqrt{\lambda/B},\qquad 0\le i<N.
\end{equation}

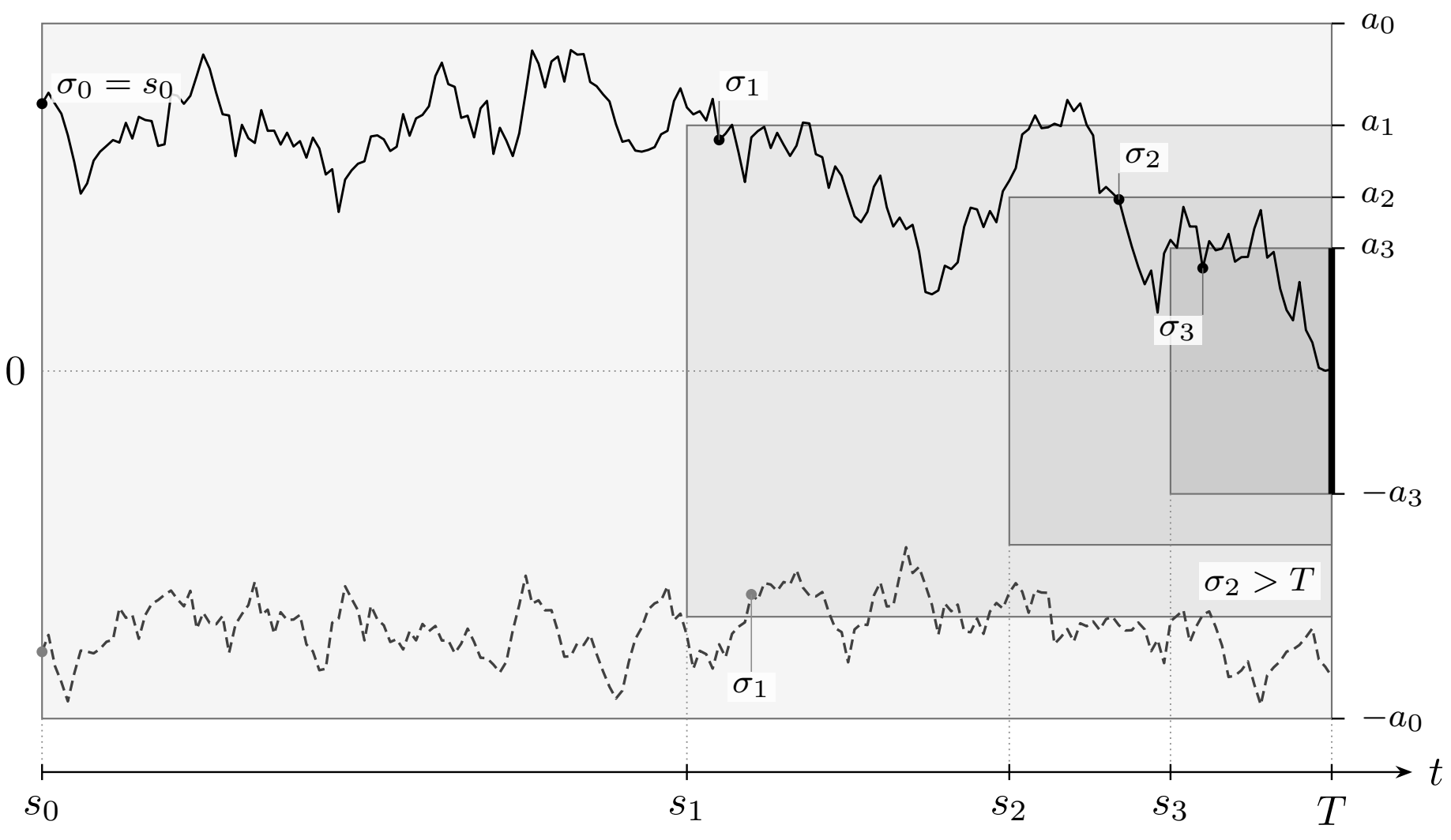
\begin{figure}[htbp]
\centering
\begin{tikzpicture}[
  x=0.9cm,y=1.12cm,font=\small,>=Stealth,
  axis/.style={line width=0.45pt,-{Stealth[length=4pt]}},
  guide/.style={black!45,line width=0.35pt,densely dotted},
  hit/.style={black,line width=0.55pt,line join=round},
  miss/.style={black!75,line width=0.62pt,dash pattern=on 2.5pt off 1.4pt,line join=round},
  lab/.style={font=\scriptsize,fill=white,fill opacity=0.85,text opacity=1,inner sep=1.2pt}
]
  \fill[black!4]  (0,-2.600) rectangle (12,2.600);
  \fill[black!9]  (6,-1.838) rectangle (12,1.838);
  \fill[black!14] (9,-1.300) rectangle (12,1.300);
  \fill[black!20] (10.5,-0.919) rectangle (12,0.919);
  \draw[black!55,line width=0.4pt] (0,-2.600) rectangle (12,2.600);
  \draw[black!55,line width=0.4pt] (6,-1.838) rectangle (12,1.838);
  \draw[black!55,line width=0.4pt] (9,-1.300) rectangle (12,1.300);
  \draw[black!55,line width=0.4pt] (10.5,-0.919) rectangle (12,0.919);
  \draw[guide] (0,0) -- (12,0);
  \node[left] at (0,0) {$0$};

  \draw[axis] (0,-3.0) -- (12.75,-3.0) node[right] {$t$};
  \foreach \x/\lab in {0/$s_0$,6/$s_1$,9/$s_2$,10.5/$s_3$,12/$T$}{
    \draw[line width=0.45pt] (\x,-2.94) -- (\x,-3.06);
    \node[below] at (\x,-3.06) {\lab};
  }
  \foreach \x/\y in {0/2.600,6/1.838,9/1.300,10.5/0.919,12/2.600}{\draw[guide] (\x,-\y) -- (\x,-2.94);}

  \foreach \y/\lab in {2.600/$a_0$,1.838/$a_1$,1.300/$a_2$,0.919/$a_3$,-0.919/$-a_3$,-2.600/$-a_0$}{
    \draw[line width=0.45pt] (12,\y) -- (12.12,\y);
    \node[right,font=\scriptsize] at (12.12,\y) {\lab};
  }

  \draw[hit] plot coordinates {
    (0.00,2.000) (0.06,2.082) (0.12,2.000) (0.18,1.925) (0.24,1.764) (0.30,1.563)
    (0.36,1.327) (0.42,1.404) (0.48,1.574) (0.54,1.641) (0.60,1.684) (0.66,1.728)
    (0.72,1.709) (0.78,1.857) (0.84,1.739) (0.90,1.902) (0.96,1.878) (1.02,1.870)
    (1.08,1.684) (1.14,1.696) (1.20,2.071) (1.26,2.060) (1.32,1.999) (1.38,2.059)
    (1.44,2.210) (1.50,2.368) (1.56,2.261) (1.62,2.082) (1.68,1.921) (1.74,1.912)
    (1.80,1.607) (1.86,1.842) (1.92,1.741) (1.98,1.706) (2.04,1.952) (2.10,1.797)
    (2.16,1.798) (2.22,1.694) (2.28,1.783) (2.34,1.680) (2.40,1.720) (2.46,1.596)
    (2.52,1.746) (2.58,1.665) (2.64,1.470) (2.70,1.510) (2.76,1.190) (2.82,1.432)
    (2.88,1.501) (2.94,1.551) (3.00,1.569) (3.06,1.755) (3.12,1.763) (3.18,1.733)
    (3.24,1.646) (3.30,1.678) (3.36,1.922) (3.42,1.760) (3.48,1.890) (3.54,1.916)
    (3.60,1.981) (3.66,2.207) (3.72,2.307) (3.78,2.146) (3.84,2.126) (3.90,1.894)
    (3.96,1.910) (4.02,1.748) (4.08,1.965) (4.14,2.021) (4.20,1.624) (4.26,1.820)
    (4.32,1.724) (4.38,1.608) (4.44,1.778) (4.50,2.082) (4.56,2.398) (4.62,2.299)
    (4.68,2.122) (4.74,2.316) (4.80,2.353) (4.86,2.165) (4.92,2.401) (4.98,2.367)
    (5.04,2.371) (5.10,2.162) (5.16,2.131) (5.22,2.070) (5.28,2.017) (5.34,1.847)
    (5.40,1.715) (5.46,1.727) (5.52,1.649) (5.58,1.641) (5.64,1.654) (5.70,1.675)
    (5.76,1.771) (5.82,1.797) (5.88,2.020) (5.94,2.115) (6.00,1.972) (6.06,1.920)
    (6.12,1.944) (6.18,1.874) (6.24,2.036) (6.30,1.729) (6.36,1.776) (6.42,1.842)
    (6.48,1.638) (6.54,1.414) (6.60,1.748) (6.66,1.795) (6.72,1.827) (6.78,1.668)
    (6.84,1.782) (6.90,1.693) (6.96,1.609) (7.02,1.686) (7.08,1.859) (7.14,1.852)
    (7.20,1.622) (7.26,1.600) (7.32,1.370) (7.38,1.531) (7.44,1.460) (7.50,1.304)
    (7.56,1.158) (7.62,1.113) (7.68,1.191) (7.74,1.379) (7.80,1.462) (7.86,1.225)
    (7.92,1.081) (7.98,1.148) (8.04,1.061) (8.10,1.094) (8.16,0.896) (8.22,0.591)
    (8.28,0.574) (8.34,0.603) (8.40,0.788) (8.46,0.763) (8.52,0.813) (8.58,1.079)
    (8.64,1.222) (8.70,1.209) (8.76,1.077) (8.82,1.197) (8.88,1.113) (8.94,1.345)
    (9.00,1.424) (9.06,1.518) (9.12,1.770) (9.18,1.809) (9.24,1.912) (9.30,1.817)
    (9.36,1.823) (9.42,1.850) (9.48,1.832) (9.54,2.028) (9.60,1.944) (9.66,2.002)
    (9.72,1.841) (9.78,1.763) (9.84,1.332) (9.90,1.377) (9.96,1.331) (10.02,1.283)
    (10.08,1.101) (10.14,0.931) (10.20,0.779) (10.26,0.649) (10.32,0.752) (10.38,0.437)
    (10.44,0.880) (10.50,0.981) (10.56,0.922) (10.62,1.228) (10.68,1.081) (10.74,1.081)
    (10.80,0.770) (10.86,0.972) (10.92,0.903) (10.98,0.916) (11.04,1.026) (11.10,0.818)
    (11.16,0.852) (11.22,0.854) (11.28,1.066) (11.34,1.204) (11.40,0.848) (11.46,0.891)
    (11.52,0.617) (11.58,0.456) (11.64,0.379) (11.70,0.666) (11.76,0.307) (11.82,0.213)
    (11.88,0.024) (11.94,0.002) (12.00,0.014)
  };
  \draw[miss] plot coordinates {
    (0.00,-2.100) (0.06,-1.974) (0.12,-2.206) (0.18,-2.327) (0.24,-2.471) (0.30,-2.262)
    (0.36,-2.091) (0.42,-2.099) (0.48,-2.112) (0.54,-2.078) (0.60,-2.025) (0.66,-2.010)
    (0.72,-1.772) (0.78,-1.846) (0.84,-1.825) (0.90,-2.003) (0.96,-1.826) (1.02,-1.742)
    (1.08,-1.713) (1.14,-1.674) (1.20,-1.642) (1.26,-1.708) (1.32,-1.760) (1.38,-1.644)
    (1.44,-1.928) (1.50,-1.808) (1.56,-1.885) (1.62,-1.894) (1.68,-1.832) (1.74,-2.109)
    (1.80,-1.896) (1.86,-1.818) (1.92,-1.751) (1.98,-1.576) (2.04,-1.826) (2.10,-1.785)
    (2.16,-1.966) (2.22,-1.804) (2.28,-1.862) (2.34,-1.859) (2.40,-1.824) (2.46,-2.045)
    (2.52,-2.097) (2.58,-2.238) (2.64,-2.227) (2.70,-1.882) (2.76,-1.851) (2.82,-1.610)
    (2.88,-1.703) (2.94,-1.791) (3.00,-2.015) (3.06,-1.755) (3.12,-1.858) (3.18,-1.902)
    (3.24,-2.027) (3.30,-2.006) (3.36,-2.083) (3.42,-1.948) (3.48,-2.015) (3.54,-1.892)
    (3.60,-1.946) (3.66,-1.907) (3.72,-2.016) (3.78,-2.015) (3.84,-2.111) (3.90,-2.040)
    (3.96,-1.927) (4.02,-2.026) (4.08,-2.145) (4.14,-2.148) (4.20,-2.196) (4.26,-2.254)
    (4.32,-2.175) (4.38,-1.941) (4.44,-1.751) (4.50,-1.531) (4.56,-1.740) (4.62,-1.715)
    (4.68,-1.790) (4.74,-1.789) (4.80,-1.947) (4.86,-2.137) (4.92,-2.133) (4.98,-2.037)
    (5.04,-2.049) (5.10,-1.974) (5.16,-2.117) (5.22,-2.246) (5.28,-2.362) (5.34,-2.450)
    (5.40,-2.390) (5.46,-2.174) (5.52,-2.005) (5.58,-1.919) (5.64,-1.784) (5.70,-1.738)
    (5.76,-1.717) (5.82,-1.609) (5.88,-1.863) (5.94,-1.815) (6.00,-1.978) (6.06,-2.227)
    (6.12,-2.093) (6.18,-2.118) (6.24,-2.226) (6.30,-2.044) (6.36,-2.141) (6.42,-1.966)
    (6.48,-1.913) (6.54,-1.880) (6.60,-1.671) (6.66,-1.714) (6.72,-1.586) (6.78,-1.597)
    (6.84,-1.645) (6.90,-1.578) (6.96,-1.596) (7.02,-1.489) (7.08,-1.620) (7.14,-1.679)
    (7.20,-1.686) (7.26,-1.653) (7.32,-1.804) (7.38,-1.923) (7.44,-1.953) (7.50,-2.179)
    (7.56,-1.933) (7.62,-1.892) (7.68,-1.900) (7.74,-1.679) (7.80,-1.581) (7.86,-1.761)
    (7.92,-1.756) (7.98,-1.530) (8.04,-1.318) (8.10,-1.511) (8.16,-1.466) (8.22,-1.603)
    (8.28,-1.748) (8.34,-1.972) (8.40,-1.739) (8.46,-1.797) (8.52,-1.747) (8.58,-1.946)
    (8.64,-1.954) (8.70,-1.851) (8.76,-1.966) (8.82,-1.799) (8.88,-1.735) (8.94,-1.776)
    (9.00,-1.658) (9.06,-1.587) (9.12,-1.650) (9.18,-1.808) (9.24,-1.640) (9.30,-1.657)
    (9.36,-1.659) (9.42,-2.042) (9.48,-1.991) (9.54,-1.930) (9.60,-2.024) (9.66,-1.887)
    (9.72,-1.909) (9.78,-1.873) (9.84,-1.936) (9.90,-1.856) (9.96,-1.838) (10.02,-1.902)
    (10.08,-1.941) (10.14,-1.939) (10.20,-1.881) (10.26,-1.933) (10.32,-2.096) (10.38,-2.010)
    (10.44,-2.185) (10.50,-1.872) (10.56,-1.831) (10.62,-1.788) (10.68,-2.028) (10.74,-1.910)
    (10.80,-1.824) (10.86,-1.798) (10.92,-1.914) (10.98,-2.055) (11.04,-2.289) (11.10,-2.278)
    (11.16,-2.239) (11.22,-2.172) (11.28,-2.346) (11.34,-2.497) (11.40,-2.275) (11.46,-2.217)
    (11.52,-2.172) (11.58,-2.099) (11.64,-2.083) (11.70,-2.049) (11.76,-1.991) (11.82,-1.925)
    (11.88,-2.159) (11.94,-2.212) (12.00,-2.285)
  };

  \draw[line width=1.6pt] (12,-0.919) -- (12,0.919);

  \fill (0,2.000) circle[radius=1.3pt];
  \node[lab,anchor=south west] at (0.08,2.000) {$\sigma_0=s_0$};
  \fill (6.30,1.729) circle[radius=1.3pt];
  \node[lab,anchor=south west] at (6.30,2.02) {$\sigma_1$};
  \draw[black!60,line width=0.3pt] (6.30,1.729) -- (6.30,2.02);
  \fill (10.02,1.283) circle[radius=1.3pt];
  \node[lab,anchor=south west] at (10.02,1.48) {$\sigma_2$};
  \draw[black!60,line width=0.3pt] (10.02,1.283) -- (10.02,1.48);
  \fill (10.80,0.770) circle[radius=1.3pt];
  \node[lab,anchor=north east] at (10.80,0.42) {$\sigma_3$};
  \draw[black!60,line width=0.3pt] (10.80,0.770) -- (10.80,0.42);
  \fill[black!50] (0,-2.100) circle[radius=1.3pt];
  \fill[black!50] (6.60,-1.671) circle[radius=1.3pt];
  \node[lab,anchor=north,text=black] at (6.60,-2.25) {$\sigma_1$};
  \draw[black!45,line width=0.3pt] (6.60,-1.671) -- (6.60,-2.25);
  \node[lab,anchor=south east,text=black] at (11.9,-1.72) {$\sigma_2>T$};
\end{tikzpicture}
\caption{The nested rectangles $[s_i,T]\times[-a_i,a_i]$,
$a_i=\gamma\sqrt{\lambda k_i}$, in the proof of
Proposition~\ref{lem:rectangle}, for $N=3$. A path ending on the thick segment at time $T$, where
$|M_T|\le a_N$, must enter every rectangle. The solid path shows the successive entrance
times $\sigma_0,\ldots,\sigma_3$. The dashed path enters the first two
rectangles and misses the third, so that $\sigma_2>T$.}
\label{fig:nested-rectangles}
\end{figure}

We first estimate the probability of avoiding the inner spatial interval
when the martingale is within one jump of the boundary of the outer
interval. The bound $b\le (\varepsilon/\gamma)a_N\le a_N/10$ and
the inequality $a_{i+1}\le a_i/\sqrt{2}$ give
\begin{equation}\label{eq:rectangle-gap}
 a_i-b-a_{i+1}\ge a_i/10,\qquad 0\le i<N.
\end{equation}
Fix $i<N$ and a stopping time $s\in[s_i,T]$. Work conditionally on
$\FF_s$ on the event $a_i-b\le |M_s|\le a_i+b$. For $t\ge s$ set
$X_t=\operatorname{sgn}(M_s)M_t$, so that $X_s=m:=|M_s|$ and $X$ has
the same conditional increment variances as $M$ after $s$. Put
$a=a_{i+1}$ and $K=a+4\sqrt{Bk_i}$. Then $a<m<K$. Let $\xi$ be the
first exit of $X$ from $(a,K)$ during $[s,T]$, capped at $T$, and
write $G=\{X_\xi>a\}$. On $G$ the path either exits at the upper end or
remains in $(a,K)$ until time $T$. Since $X_\xi\le a$ on $G^c$ and
$X_\xi\le K+b$ on $G$, optional stopping gives
\[
 m=\EE[X_\xi\mid\FF_s]
 \le a+(K+b-a)\PP(G\mid\FF_s).
\]
On $G\cap\{\xi=T\}$ the path has stayed above $a$ through time $T$.
On $G\cap\{\xi<T\}$ we have $X_\xi\ge K$, so a subsequent visit to
$(-\infty,a]$ requires a displacement of at least $K-a$. Doob's maximal
inequality for the squared displacement, followed by
\eqref{eq:rectangle-bracket} at $\xi$, bounds its conditional
probability given $\FF_\xi$ by
\[
 \frac{\EE[(X_T-X_\xi)^2\mid\FF_\xi]}{(K-a)^2}
 \le \frac{Bk_i}{(K-a)^2}=\frac1{16}.
\]
Thus, on $G$, with conditional probability at least $15/16$ the path
stays above $a$ from $s$ through $T$. Taking conditional expectations
at $s$ yields
\begin{equation}\label{eq:rectangle-escape}
 \PP\bigl(|M_t|>a_{i+1}\text{ for all }s\le t\le T\mid\FF_s\bigr)
 \ge \frac{15}{16}\frac{m-a}{K+b-a}\ge c_3\sqrt{\lambda/B}.
\end{equation}
Indeed, $m-a\ge a_i/10$ by \eqref{eq:rectangle-gap}, whereas
$K+b-a\le5\sqrt{Bk_i}$. In particular $c_3=\gamma/100$ suffices.

We now apply this estimate at the first entrance into the outer
rectangle. Suppose first that $\sigma_i=s_i$. The martingale then starts
the block in $[-a_i,a_i]$. Let $\tau$ be the first time $t\ge s_i$ with
$|M_t|>a_i$, capped at $s_{i+1}$. The event
$E=\{|M_\tau|>a_i\}\in\FF_\tau$ records an exit by $s_{i+1}$.
If $E^c\cap A_i$ occurs, then $\tau=s_{i+1}$ and
$V_\tau-V_{s_i}\ge\lambda k_i$. Also $|M_\tau|\le a_i+b$, even on an
exit. Optional stopping at $\tau$ gives
\[
 \begin{aligned}
 \lambda k_i\PP(E^c\cap A_i\mid\FF_{s_i})
 &\le \EE[V_\tau-V_{s_i}\mid\FF_{s_i}]
   = \EE[(M_\tau-M_{s_i})^2\mid\FF_{s_i}]\\
 &\le (2a_i+b)^2\le(2\gamma+\varepsilon)^2\lambda k_i
   <\frac12\lambda k_i .
 \end{aligned}
\]
Thus $E$ and $E^c\cap A_i^c$ have combined conditional probability
at least $1/2$. On $E$ we have $a_i<|M_\tau|\le a_i+b$. Since
$\tau\le s_{i+1}$, avoiding $[-a_{i+1},a_{i+1}]$ from $\tau$ through
$T$ implies $\sigma_{i+1}>T$. Hence \eqref{eq:rectangle-escape} at
$\tau$ gives $\PP(\sigma_{i+1}>T\mid\FF_\tau)\ge p$ on $E$, where
$p=c_3\sqrt{\lambda/B}\le1$. The events $E\cap\{\sigma_{i+1}>T\}$ and
$E^c\cap A_i^c$ are disjoint. Conditioning first at $\tau$ in the
first event gives
\[
 \PP(\{\sigma_{i+1}>T\}\cup A_i^c\mid\FF_{s_i})
 \ge p\,\PP(E\mid\FF_{s_i})+\PP(E^c\cap A_i^c\mid\FF_{s_i})
 \ge p/2.
\]
This proves \eqref{eq:rectangle-miss} with $c_1=c_3/2$.

If $s_i<\sigma_i\le T$, the martingale enters from outside
$[-a_i,a_i]$. The jump bound therefore places $|M_{\sigma_i}|$ between
$a_i-b$ and $a_i$. No inner rectangle has been entered before
$\sigma_i$, since it is contained in the outer rectangle. Applying
\eqref{eq:rectangle-escape} at the bounded stopping time
$s=\sigma_i\wedge T$ on this event gives \eqref{eq:rectangle-miss}
also in this case.

Finally put $E_i=\{\sigma_i\le T\}\cap\bigcap_{j<i}A_j$. At time
$\sigma_i\wedge T$ we know whether the $i$th rectangle was entered
and whether each of the earlier bounds $A_j$ held, since
$s_{j+1}\le s_i$ for $j<i$. As $E_{i+1}\subseteq E_i$, we may therefore
condition at $\sigma_i\wedge T$ in \eqref{eq:rectangle-miss} to obtain
\[
 \PP(E_{i+1})\le(1-c_1\sqrt{\lambda/B})\PP(E_i).
\]
Take $\eta=\gamma$. On $\{|M_T|\le\eta\sqrt{\lambda k_N}\}$ the
martingale has entered the last rectangle. Intersecting with all the
events $A_i$ thus gives a subset of $E_N$, and iteration proves
\eqref{eq:abstract-smallball}.
\end{proof}

In the application to the height, the jump bound will hold on an event
determined by the planar path, which is revealed at time zero. We shall
also have an upper bound for the total bracket on each terminal window
outside an event of small probability. The following corollary allows
both restrictions to appear in the probability being estimated.

\begin{corollary}\label{cor:localized}
Let $M$ be a square-integrable martingale. Take the times, $\lambda$,
$B$ and absolute constants of Proposition \ref{lem:rectangle}. Define
$A_i$ by \eqref{eq:abstract-lower-events} and put
$G_i=\{V_T-V_{s_i}\le Bk_i\}$. If \eqref{eq:abstract-jumps} holds
on an event $E\in\FF_0$, then for every $z\in\RR$,
\[
 \PP\!\left(|M_T-z|\le\eta\sqrt{\lambda k_N},\
       E\cap\bigcap_{i<N}(A_i\cap G_i)\right)
       \le e^{-cN\sqrt{\lambda/B}} .
\]
\end{corollary}

\begin{proof}
It suffices to consider $\PP(E)>0$ and to work under
$\PP(\,\cdot\mid E)$. Since $E\in\FF_0$, the martingale property and
the conditional increment variances agree with the original ones on
$E$. For each terminal window that has already started, inspect its
bracket before taking the next step. Let $\nu$ be the first $r<T$ for
which $V_{r+1}-V_{s_i}>Bk_i$ for some $i<N$ with $s_i\le r$, and put
$\nu=T$ if there is no such $r$. This is a stopping time because
$V_{r+1}$ is $\FF_r$-measurable.

The stopped martingale has bracket $V_{t\wedge\nu}$. If $\nu\le s_i$,
this bracket is constant throughout $[s_i,T]$. If $\nu>s_i$, the
test at $r=\nu-1$ gave $V_\nu-V_{s_i}\le Bk_i$. Thus in both cases
the total increase of the stopped bracket on $[s_i,T]$ is at most
$Bk_i$ pathwise. Its increase after any $t\in[s_i,T]$ obeys the same
bound by monotonicity, which gives \eqref{eq:rectangle-bracket} for
the stopped martingale. The jump bound holds under the conditional law.
On $\bigcap_{i<N}G_i$ none of the tests can fail, so $\nu=T$ and both
the terminal value and every event $A_i$ agree with their stopped
counterparts. Applying Proposition \ref{lem:rectangle} and multiplying
by $\PP(E)\le1$ gives the claimed bound.
\end{proof}

\section{Transience and range}\label{sec:transience}

\subsection{Transience}

Fix $\theta\in(0,1]$. Excursion $n$ lies in column $U_n$ and ends at
height $H_n$. To bound the probability that it visits a given vertex,
we first estimate the chance that $H_n$ lies near a prescribed height
given $U_n=u$. As in Section \ref{sec:bridge}, we multiply this
conditional probability by $(n+1)\PP(U_n=u)$. This gives an estimate
uniform in the choice of vertex before we sum over $n$.

\begin{proposition}\label{prop:height-smallball}
For every $0<\beta<1/2$, $A>0$ and $n\ge1$,
\begin{equation}\label{eq:height-smallball}
 \sup_{u\in\ZZ^2,\ z\in\RR}
 (n+1)\PP(U_n=u,\ |H_n-z|\le n^\beta)
              \le C_{\theta,\beta,A} L_n^{-A}.
\end{equation}
At any fixed $u$ this gives
$\sup_z Q_n^u(|H_n-z|\le n^\beta)\le C_{\theta,\beta,A,u}L_n^{-A}$
at every admissible time.
\end{proposition}

At the endpoint $u=(n,0)$, all $n$ planar increments equal $(1,0)$.
Every arrival is fresh, so $Q_n^{(n,0)}(H_n=0)=1$. The factor
$(n+1)\PP(U_n=u)$ in \eqref{eq:height-smallball} is then only
$(n+1)(\theta/4)^n$. This example explains the dependence on $u$ in
the conditional estimate.

\begin{proof}[Proof of Proposition \ref{prop:height-smallball}]
Fix an admissible $(n,u)$ and write $Q=Q_n^u$, $\mu=\Pi_n^u$ and $L=L_n$.
To apply Corollary \ref{cor:localized}, we need a lower bound for the
bracket on each of a family of disjoint blocks, an upper bound on the
whole window from the beginning of each block to $n$, and a bound on
the jumps. 

In Section \ref{sec:column} we revealed the planar increments successively
to estimate the contribution of first departures to the clock. Here we
use $\mathcal H_j=\mathcal D_j\vee\sigma(U)$, so the whole planar path is
known at time zero. By Lemma \ref{lem:height-martingale}(c), $H$ is still
a martingale under $Q$. Its bracket $V$ has increments
\begin{equation}\label{eq:height-variance}
 v_j:=V_{j+1}-V_j
    =\EE_Q[D_{j+1}\mid\mathcal H_j]\ge\ind_{\{D_{j+1}>0\}} .
\end{equation}
Indeed, $D_{j+1}>0$ exactly when $I_{j+1}$ is nonempty, and this event
belongs to $\mathcal A_{j+1}\subseteq\mathcal H_j$.

We next choose the dyadic blocks and terminal windows to which
Corollary \ref{cor:localized} will be applied. Set $k_0=\lfloor n/2\rfloor$, $k_i=\lfloor2^{-i}k_0\rfloor$ and
$s_i=n-k_i$. Let $J$ be the last index with $k_J\ge n^{2\beta}$.
For all large $n$ we have $J\ge1$ and
$J=(1-2\beta)\log n/\log2+O(1)$. For $i<J$, the $i$th block consists
of excursions $s_i+1,\ldots,s_{i+1}$, while the corresponding terminal
window consists of excursions $s_i+1,\ldots,n$
(Figure \ref{fig:timeline}). The blocks use only the last $k_0$ planar
increments. The final length $k_J$ gives a spatial scale $\sqrt{k_J}$
of order $n^\beta$ in Proposition \ref{lem:rectangle}.

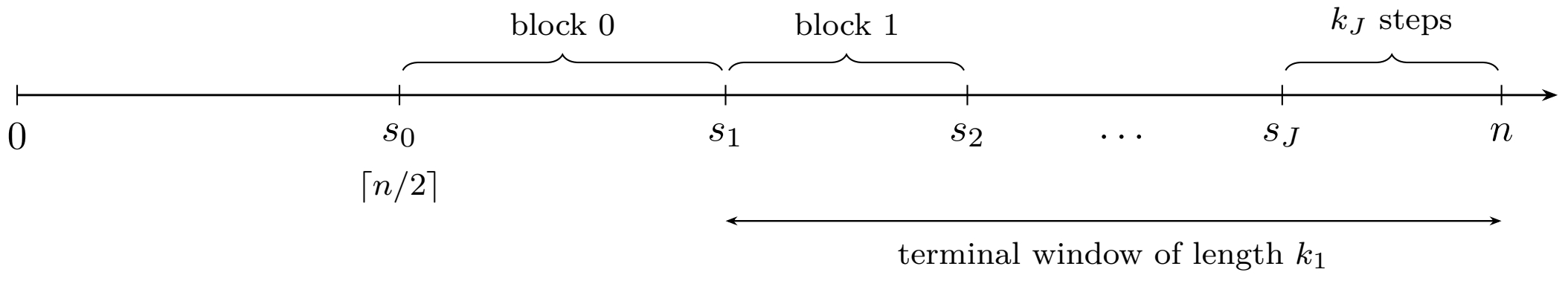
\begin{figure}[!htbp]
\centering
\begin{tikzpicture}[x=1cm,y=1cm,>=Stealth,font=\small]
 \draw[-{Stealth[length=4pt]},line width=.55pt] (0,0) -- (13.7,0);
 \foreach \x/\lab in {0/$0$,3.4/$s_0$,6.3/$s_1$,8.45/$s_2$,11.25/$s_J$,13.2/$n$}{
   \draw (\x,-.09) -- (\x,.09);
   \node[below] at (\x,-.12) {\lab};
 }
 \node[below,font=\scriptsize] at (3.4,-.54) {$\lceil n/2\rceil$};
 \node at (9.85,-.38) {$\cdots$};
 \draw[decorate,decoration={brace,amplitude=4pt}] (3.43,.22) -- (6.27,.22);
 \draw[decorate,decoration={brace,amplitude=4pt}] (6.33,.22) -- (8.42,.22);
 \node[above,font=\scriptsize] at (4.85,.40) {block $0$};
 \node[above,font=\scriptsize] at (7.38,.40) {block $1$};
 \draw[decorate,decoration={brace,amplitude=4pt}] (11.28,.22) -- (13.17,.22);
 \node[above,font=\scriptsize] at (12.22,.40) {$k_J$ steps};
 \draw[{Stealth[length=3pt]}-{Stealth[length=3pt]},line width=.45pt] (6.3,-1.12) -- (13.2,-1.12);
 \node[below,font=\scriptsize] at (9.75,-1.17) {terminal window of length $k_1$};
\end{tikzpicture}
\caption{Disjoint pairs give the lower bound for $V_{s_{i+1}}-V_{s_i}$. Moment bounds for $T_n-T_{s_i}$ yield moment bounds for $V_n-V_{s_i}$ through \eqref{eq:predictable-projection}. The last length $k_J$ is at least
$n^{2\beta}$.}
\label{fig:timeline}
\end{figure}

For the lower bounds, we use deterministic patterns in the planar
increments which force vertical motion during at least one excursion of
a pair. Pair the successive planar increments in each block, dropping the last increment if necessary. Mark $(Y_{j+1},Y_{j+2})$ if $Y_{j+1}=0$, or if
$Y_{j+1}\ne0$ and $Y_{j+2}=-Y_{j+1}$. Each marked pair forces at least
one vertical step. A zero increment gives $D_{j+1}\ge1$. For a
reversal, if $D_{j+1}=0$, the first step discovers
$W_{j+1}=(U_{j+1},H_j)$ and the second returns to the consumed $W_j$,
forcing $D_{j+2}\ge1$. If $D_{j+1}>0$, the first excursion already
contributes. Thus $v_j+v_{j+1}\ge1$ for each marked pair by
\eqref{eq:height-variance}. Figure \ref{fig:pair-reversal} also shows
why an individual excursion can have zero variance.

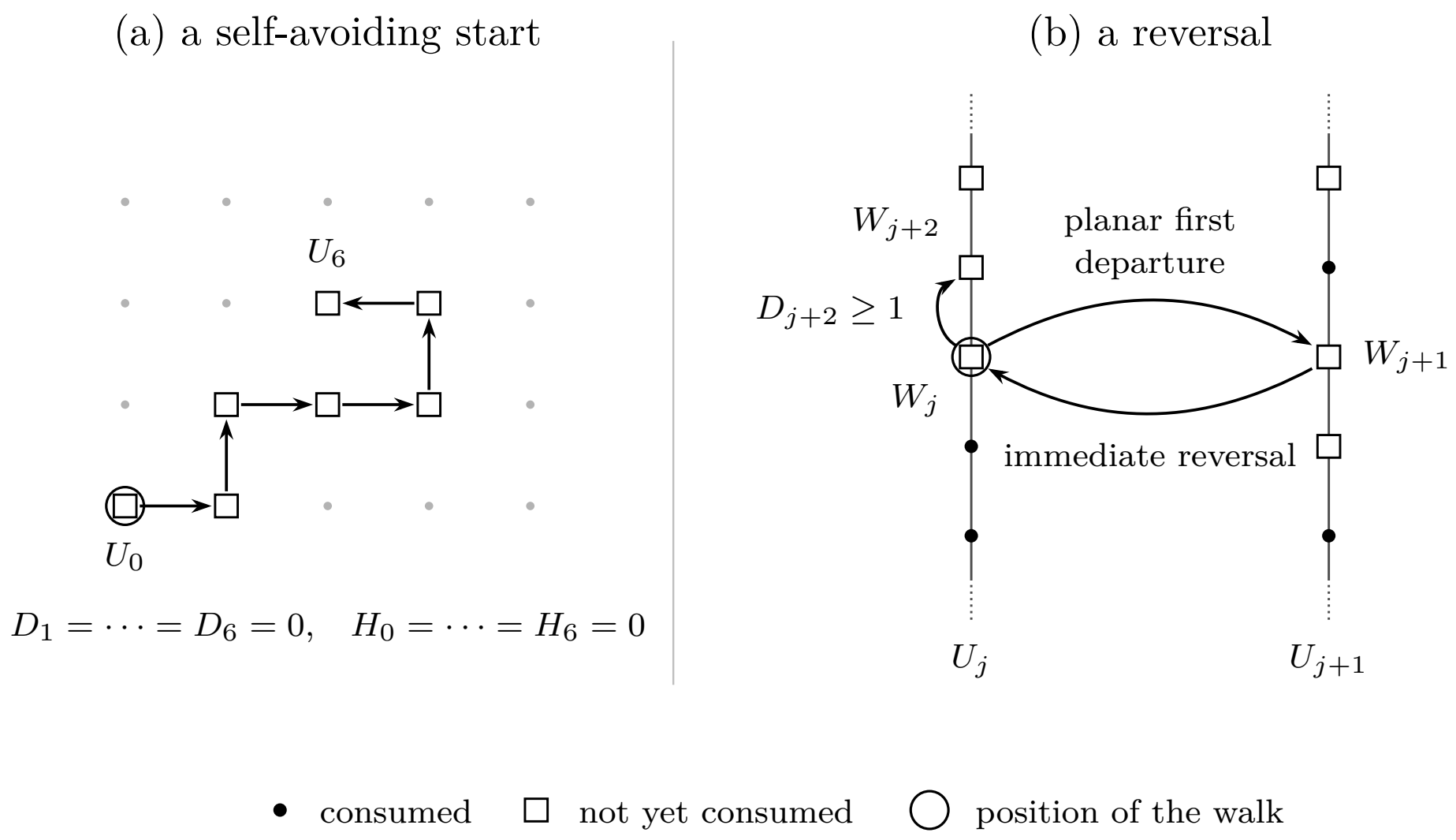
\begin{figure}[!tbp]
\centering
\begin{tikzpicture}[
    >=Stealth,
    font=\small,
    consumed/.style={circle,fill=black,inner sep=0pt,minimum size=3.2pt},
    unvisited/.style={rectangle,draw=black,fill=white,line width=0.55pt,
                     inner sep=0pt,minimum size=5.4pt},
    walker/.style={circle,draw=black,line width=0.6pt,inner sep=0pt,
                   minimum size=9pt},
    move/.style={->,line width=0.7pt,
                 >={Stealth[length=1.7mm,width=1.15mm]}},
    lab/.style={font=\scriptsize}
]
\begin{scope}[x=0.85cm,y=0.85cm]
\node at (2.0,4.65) {(a) a self-avoiding start};
\foreach \x in {0,...,4}{\foreach \y in {0,...,3}{
  \fill[black!30] (\x,\y) circle[radius=1pt];}}
\foreach \x/\y in {0/0,1/0,1/1,2/1,3/1,3/2,2/2}{\node[unvisited] at (\x,\y) {};}
\node[walker] at (0,0) {};
\foreach \a/\b/\c/\d in {0/0/1/0,1/0/1/1,1/1/2/1,2/1/3/1,3/1/3/2,3/2/2/2}{
  \draw[move,shorten <=3.5pt,shorten >=3.5pt] (\a,\b) -- (\c,\d);
}
\node[lab,anchor=north] at (0,-0.2) {$U_0$};
\node[lab,anchor=south] at (2,2.2) {$U_6$};
\node[lab,align=center] at (2.0,-1.2)
  {$D_1=\cdots=D_6=0$,\quad $H_0=\cdots=H_6=0$};
\end{scope}

\draw[black!25,line width=0.45pt] (4.6,-1.5) -- (4.6,3.9);

\begin{scope}[shift={(7.1,1.25)},x=1cm,y=0.75cm]
\node at (1.5,3.6) {(b) a reversal};
\foreach \x in {0,3}{
  \draw[black!70,line width=0.55pt] (\x,-2.5) -- (\x,2.5);
  \draw[densely dotted,black!70,line width=0.55pt] (\x,2.55) -- (\x,2.95);
  \draw[densely dotted,black!70,line width=0.55pt] (\x,-2.55) -- (\x,-2.95);
}
\foreach \h in {-2,-1}{\node[consumed] at (0,\h) {};}
\foreach \h in {0,1,2}{\node[unvisited] at (0,\h) {};}
\node[walker] at (0,0) {};
\foreach \h in {-2,1}{\node[consumed] at (3,\h) {};}
\foreach \h in {-1,0,2}{\node[unvisited] at (3,\h) {};}
\draw[move] (0.14,0.13) to[bend left=28] (2.86,0.13);
\draw[move] (2.86,-0.13) to[bend left=28] (0.14,-0.13);
\node[lab,align=center] at (1.5,1.27) {planar first\\departure};
\node[lab,align=center,anchor=north] at (1.5,-0.78)
  {immediate reversal};
\draw[move] (-0.13,0.13) to[bend left=60] (-0.13,0.87);
\node[lab,anchor=east] at (-0.42,0.5) {$D_{j+2}\ge1$};
\node[lab,anchor=north east] at (-0.12,-0.12) {$W_{j}$};
\node[lab,anchor=south east] at (-0.12,1.12) {$W_{j+2}$};
\node[lab,anchor=west] at (3.15,0) {$W_{j+1}$};
\node[lab] at (0,-3.4) {$U_{j}$};
\node[lab] at (3,-3.4) {$U_{j+1}$};
\end{scope}

\node[consumed] at (1.3,-2.55) {};
\node[anchor=west,font=\scriptsize,text height=1.5ex,text depth=.25ex] at (1.5,-2.55) {consumed};
\node[unvisited] at (3.45,-2.55) {};
\node[anchor=west,font=\scriptsize,text height=1.5ex,text depth=.25ex] at (3.67,-2.55) {not yet consumed};
\node[walker] at (6.75,-2.55) {};
\node[anchor=west,font=\scriptsize,text height=1.5ex,text depth=.25ex] at (7.0,-2.55) {position of the walk};
\end{tikzpicture}
\caption{(a) A self-avoiding start of the planar path for $\theta=1$
produces zero height increments. Each square shows a vertex at the
time of its discovery. (b) A planar first departure from
$W_j$ reaches a fresh vertex $W_{j+1}$. An immediate reversal returns
to the consumed $W_j$ and forces $D_{j+2}\ge1$. The vertical exit shown
is one possible realization.}
\label{fig:pair-reversal}
\end{figure}

Let $N_i$ count the marks in block $i$. That block has
$k_i-k_{i+1}\ge k_i/2$ increments, hence at least $k_i/8$ disjoint
pairs when $n$ is large. Under $\PP$ their marks are independent, each
with probability $1-\theta+\theta^2/4=(1-\theta/2)^2\ge1/4$.
A binomial lower-tail bound gives $\PP(N_i<k_i/64)\le Ce^{-ck_i}$.
The marks depend only on the last $k_0$ increments, so the bounded
density in \eqref{eq:block-density} gives the same bound, up to a
constant, for their exceptional $\mu$-mass. Since
$V_{s_{i+1}}-V_{s_i}\ge N_i$ for every realization of the vertical
signs, a union bound gives, with $\lambda=1/64$,
\begin{equation}
 \mu\!\left(\exists\,i<J:\ V_{s_{i+1}}-V_{s_i}<\lambda k_i\right)
 \le C J e^{-c n^{2\beta}} .
\end{equation}

For the upper bounds, the clock is the sum of the actual durations,
whereas the bracket is the sum of their conditional means. We therefore
use a predictable projection inequality to pass from the clock estimate
to the bracket. For $q>1$ and nonnegative adapted $X_{j+1}$, the predictable projection
inequality gives
\begin{equation}\label{eq:predictable-projection}
 \left\|\sum_j\EE_Q[X_{j+1}\mid\mathcal H_j]\right\|_{q,Q}
 \le q\left\|\sum_jX_{j+1}\right\|_{q,Q}.
\end{equation}
Indeed, for $Z\ge0$ of norm one in $L^{q/(q-1)}(Q)$, conditioning each
term on $\mathcal H_j$ gives
\[
 \EE_Q\!\left[Z\sum_j\EE_Q[X_{j+1}\mid\mathcal H_j]\right]
 =\sum_j\EE_Q\!\left[X_{j+1}\EE_Q[Z\mid\mathcal H_j]\right]
 \le\EE_Q\!\left[\Big(\max_j\EE_Q[Z\mid\mathcal H_j]\Big)
                       \sum_jX_{j+1}\right].
\]
H\"older's inequality  with exponents $q/(q-1)$ and $q$, and Doob's maximal inequality bound the last
expression by $q\|\sum_j X_{j+1}\|_{q,Q}$. Taking the supremum over $Z$
proves \eqref{eq:predictable-projection}
(see also \cite[Theorem 3.2]{BDG1972}). Apply this to the durations in
the window and multiply both norms by $w_{n,u}^{1/q}$. Since
$k_i\ge n^{2\beta}$, Proposition \ref{prop:window} yields
\begin{equation}
 \|V_n-V_{s_i}\|_{q,\mu}\le q\|T_n-T_{s_i}\|_{q,\mu}
                 \le C_{q,\beta} k_i L^{3/2}\qquad(0\le i<J).
\end{equation}
Take $B=L^{7/4}$ and let $G_i=\{V_n-V_{s_i}\le Bk_i\}$. By Markov's
inequality and $J=O(L)$,
\[
 \mu\Big(\bigcup_{i<J}G_i^c\Big)
 \le C_{q,\beta}J\left(\frac{L^{3/2}}{B}\right)^q
 \le C_{q,\beta}L^{1-q/4}.
\]
The choice of $B$ gives both a power of $L$ in this tail estimate and
$J/\sqrt B$ of order $L^{1/8}$ in the defocusing estimate.

It remains to enforce the jump condition on a planar-path event whose complement has small $\mu$-mass. Every height jump through discovery $n$ is bounded by the size of its exit interval, which is at most $\ell_n^*$. Put
$b=\varepsilon\sqrt{\lambda k_J}$ and let $E=\{\ell_n^*\le b\}$.
The event $E$ belongs to $\mathcal H_0$ since it depends only on $U$.
By \eqref{eq:local-index} with $p=1$ and $k_J\ge n^{2\beta}$, Markov's
inequality gives $\mu(E^c)\le C L^2 n^{-\beta}$.

With $A_i=\{V_{s_{i+1}}-V_{s_i}\ge\lambda k_i\}$, Corollary \ref{cor:localized} gives
\[
 Q\!\left(|H_n-z|\le\eta\sqrt{\lambda k_J},\
         E\cap\bigcap_{i<J}(A_i\cap G_i)\right)
 \le e^{-cJ\sqrt{\lambda/B}} .
\]
Multiply by $w_{n,u}\le C$ and add the exceptional masses from the
lower bracket bound, the union bound for $\bigcup_{i<J}G_i^c$, and
the estimate for $E^c$. Since $\sqrt{k_J}\ge n^\beta$ and
$J/\sqrt B\ge c_\beta L^{1/8}$, we obtain
\[
 \sup_z\mu\bigl(|H_n-z|\le\eta\sqrt{\lambda}\,n^\beta\bigr)
 \le C e^{-c_\beta L^{1/8}}+C J e^{-c n^{2\beta}}
       +C_{q,\beta}L^{1-q/4}+C L^2 n^{-\beta}.
\]
Choose $q>4(A+1)$. Each term is $O(L^{-A})$, and a fixed number of
intervals of radius $\eta\sqrt\lambda\,n^\beta$ covers any interval
of radius $n^\beta$. This proves \eqref{eq:height-smallball} for
large $n$. Increasing the constant covers the finitely many remaining
$n$, since $\mu$ has bounded mass. Inadmissible endpoints have zero
mass. For fixed $u$, division by $w_{n,u}\ge c_u$ gives the conditional
estimate at every admissible time.
\end{proof}

\begin{proof}[Proof of Theorem \ref{mainthm} and \eqref{eq:green-decay}]
Fix $x=(u,z)\in\ZZ^2\times\ZZ$ and let $E_n$ be the event that excursion
$n$ visits $x$. This requires $U_n=u$, so $E_n$ is empty if $(n,u)$ is inadmissible.
On $\{U_n=u,\ D_n>0\}$ the excursion stays in an interval of at most
$\ell_n(u)$ consumed heights until it exits at $H_n$. If $D_n=0$, it
consists only of its endpoint. In both cases
\[
 E_n\subseteq\{U_n=u,\ |H_n-z|\le\ell_n(u)\}.
\]
Write $\mu=\Pi_n^u$. By Proposition \ref{prop:height-smallball} with
$\beta=1/4$ and the fixed-site tail in \eqref{eq:local-tails}, for
every $b>0$,
\[
 \mu(E_n)\le\mu(|H_n-z|\le n^{1/4})+\mu(\ell_n(u)>n^{1/4})
            \le C_b L_n^{-b}.
\]

We also need the number of physical times at which this excursion can
visit $x$. Its duration is $D_n$, so it contains $D_n+1$ integer times
when both endpoints are counted. Given the entire planar path and the
history through $\tau_{n-1}$, Lemma \ref{lem:height-martingale}(b) bounds the conditional second moment of $D_n+1$ on
$\{U_n=u\}$ by $C(\ell_n(u)+1)^4$. Averaging and applying
Lemma \ref{lem:bridge-local} gives $\|D_n+1\|_{2,\mu}\le C L_n^2$.
For the exponent $a>0$ in \eqref{eq:green-decay}, take $b=2a+6$. Since
$E_n\subseteq\{U_n=u\}$, Cauchy--Schwarz gives
\[
 \EE[(D_n+1)\ind_{E_n}]
 = \frac{1}{n+1}\int_{E_n}(D_n+1)\,d\mu
 \le \frac{\|D_n+1\|_{2,\mu}\,\mu(E_n)^{1/2}}{n+1}
 \le \frac{C_a}{(n+1)L_n^{a+1}}.
\]

By time $\tau_n$ the walk has visited $n+1$ sites, counting $W_n$.
These sites form a connected set containing $0$, so they can contain
$x$ only if $n\ge\|x\|_1$. The excursions cover every time after zero,
and counting shared endpoints twice only increases the sum. Thus
\[
 G(x)\le\ind_{\{x=0\}}+
 \sum_{n\ge\max\{1,\|x\|_1\}}\EE[(D_n+1)\ind_{E_n}]
 \le \frac{C_a}{[\log(2+\|x\|_1)]^a}.
\]
The last inequality follows by comparison with the integral of
$1/[s(\log s)^{a+1}]$. Thus every vertex has finite expected occupation
and is visited only finitely many times almost surely. Taking the
intersection over $x\in\ZZ^3$ proves transience.
\end{proof}
\subsection{Lower bound for the range}

\begin{proof}[Proof of Proposition \ref{rangelowerbd}]
We first bound the vertical clock at discovery times. Put
$n=2^j$ and $p_j=\lceil2\log(j+2)\rceil$. Since $L_{2^j}$ grows
linearly in $j$, we have $p_j\le L_n^{1/4}$ for all sufficiently
large $j$.
By \eqref{eq:clock-high-moments} and Markov's inequality,
\[
 \PP\!\left(T_{2^j}> e C p_j^3\,2^j L_{2^j}^{3/2}\right)
       \le e^{-p_j}\le(j+2)^{-2}.
\]
The probabilities are summable in $j$, so by Borel--Cantelli the clock
obeys the indicated bound at all sufficiently large dyadic indices,
almost surely. If $2^{j-1}<n\le2^j$, then $T_n\le T_{2^j}$,
$2^j<2n$, $L_{2^j}\le C\log n$ and $p_j\le C\log\log n$ for large $n$.
This gives, almost surely for all sufficiently large $n$,
\[
 T_n\le C n(\log n)^{3/2}(\log\log n)^3.
\]
For a physical time $\tau_n<t\le\tau_{n+1}$ the first departures before
$t$ are exactly those from $W_0,\ldots,W_n$, so $R_t=n+1$. Thus
\eqref{eq:clock} and the preceding estimate give, almost surely for
all sufficiently large $t$,
\[
 t\le\tau_{n+1}\le n+1+T_{n+1}
   \le C R_t (\log R_t)^{3/2}(\log\log R_t)^3 .
\]
Since $R_t\le t$, replacing the two logarithmic factors by their
values at $t$ and dividing gives Proposition \ref{rangelowerbd}.
\end{proof}

\subsection{Lower tails of the range}\label{sec:range-lower-tail}

We now prove Proposition \ref{prop-rangeLaplacelog} directly from the
departure rules of the general walk $X=M(d,d_1,d_2)$. Write
$\FF_s=\sigma(X_0,\ldots,X_s)$ and put
$A=\{1,\ldots,d_1\}$, $B=\{d-d_2+1,\ldots,d\}$. These are the
coordinate directions of the first- and subsequent-departure laws.
Let $\ZZ^B=\{x\in\ZZ^d:x_i=0\text{ for }i\notin B\}$. We count the
current site in $N_s=R_{s+1}=|\{X_0,\ldots,X_s\}|$, so that
$N_{s+1}>N_s$ precisely when the next step discovers a new site.
We prove that there exist $r_0\ge 1$ and constants $c,C>0$ such that, for all integers $r\ge r_0$ and $t\ge0$,
\begin{equation}\label{eq:range-tail}
 \exp\!\left[-C\left(r+\frac{t}{r^{2/d_2}}\right)\right]
 \le \PP(N_{t}\le r)
 \le \min\!\left\{1,\exp\!\left(Cr-\frac{ct}{r^{2/d_2}}\right)\right\}.
\end{equation}
Here and below the constants in this subsection may depend on the
dimensions. To recover the examples in the introduction, take
$r=\lfloor t^{1/4}\rfloor$. The term $t/r^{2/d_2}$ is of order $t^{3/4}$
when $d_2=2$ and $t^{1/2}$ when $d_2=1$. Both dominate $r$.
Since $R_t=N_{t-1}$ for $t\ge1$, the same exponents hold for $R_t$.

\begin{proof}
We first prove the upper bound in \eqref{eq:range-tail}.
Put $h=1+\lceil C_0 r^{2/d_2}\rceil$. Whenever $N_{s}\le r$,
\begin{equation}\label{eq:discovery-block}
 \PP(N_{s+h}=N_{s}\mid\FF_s)\le\frac12
\end{equation}
for $C_0$ sufficiently large. Indeed, fix the history through time $s$
and write $D=\{X_0,\ldots,X_s\}$. Condition also on the first step of
the block. If $X_{s+1}\notin D$, the range has already grown. If
$X_{s+1}=y\in D$, then all vertices of $D$ have now had their first
departure, including $X_s$ itself. Until the next discovery, the walk
uses only subsequent departures and agrees with a simple random walk
$Y$ started at $y$ in $y+\ZZ^B$. Write $\PP_y$ for the law of this
simple walk. The set $\Gamma=D\cap(y+\ZZ^B)$
has at most $r$ sites. To make no discovery, $Y$ must stay in
$\Gamma$ for $h-1$ steps. In particular its last position belongs to
$\Gamma$, and the heat-kernel bound
\cite[Theorems 2.1.1 and 2.1.3]{MR2677157} gives
\[
 \PP_y(Y_{h-1}\in\Gamma)
 \le |\Gamma|\sup_{z\in y+\ZZ^B}\PP_y(Y_{h-1}=z)
 \le C r(h-1)^{-d_2/2}\le\frac12 .
\]
This choice of $C_0$ works for every history, including when $A$ and
$B$ overlap.

To iterate this bound, let
$Z_k=2^{-N_{kh}}\ind_{\{N_{kh}\le r\}}$. As long as $N_{kh}\le r$,
the next block keeps the weight unchanged if it has no discovery.
Any discovery reduces the weight by at least a factor of two, and
crossing the threshold $r$ makes it zero. Thus \eqref{eq:discovery-block}
gives $\EE[Z_{k+1}\mid\FF_{kh}]\le(\tfrac12+\tfrac12\cdot\tfrac12)Z_k
=\tfrac34 Z_k$. Once $N_{kh}>r$, both sides are zero by monotonicity.
Iteration gives $\EE Z_k\le(3/4)^k$. On $\{N_t\le r\}$ we have
$Z_k\ge2^{-r}$ for $k=\lfloor t/h\rfloor$, so
\[
 \PP(N_{t}\le r)\le 2^r\EE Z_k\le 2^r(3/4)^k
 \le\exp(Cr-ct/r^{2/d_2}).
\]

To prove the lower bound in \eqref{eq:range-tail}, let $Y$ be simple random walk on $\ZZ^B$ started at the origin.
Independently at each vertex, sample a first-departure arrow uniformly
from $\{\pm e_i:i\in A\}$. Use the increments of $Y$, independently of
these arrows, on subsequent departures of $X$.
First consider the probability that $Y$ stays in a box. For an integer
$m\ge1$, put $\Lambda_m=\{x\in\ZZ^B: |x_i|\le m\text{ for }i\in B\}$.
For $\vartheta=\pi/(2m+2)$ the function
$h_m(x)=\prod_{i\in B}\cos(\vartheta x_i)$ is positive on this box,
equals one at the origin and vanishes at every site of
$\ZZ^B\setminus\Lambda_m$ adjacent to the box.
Let $P_{\Lambda_m}$ be the transition kernel of the simple walk
killed on leaving the box. In each coordinate, averaging the values of
the cosine at the two neighboring points multiplies its value by
$\cos\vartheta$. Averaging over coordinates therefore gives
$P_{\Lambda_m}h_m=\cos(\vartheta)h_m$. Since $0<h_m\le1$ on the box,
\begin{equation}\label{eq:box-survival}
 \PP_0(Y_0,\ldots,Y_t\in\Lambda_m)
 = (P_{\Lambda_m}^t\ind_{\Lambda_m})(0)
 \ge (P_{\Lambda_m}^t h_m)(0)=(\cos\vartheta)^t
 \ge e^{-Ct/m^2}.
\end{equation}
Here $-\log\cos\vartheta\le C/m^2$ for $m\ge1$.

Fix $i\in A$ and define $\phi:\ZZ^d\to\ZZ^d$ by
$\phi(x)=x+(-1)^{x_i}e_i$. The map pairs neighboring vertices along
axis $i$ throughout $\ZZ^d$. In particular it also pairs every vertex
of $\Lambda_m$ when $i\notin B$, in which case the partner lies outside
$\ZZ^B$. Since $\phi$ is an involution,
$\mathcal P=\Lambda_m\cup\phi(\Lambda_m)$ is a union of pairs and
$|\mathcal P|\le2|\Lambda_m|$. With probability
$(2d_1)^{-|\mathcal P|}$ the first arrow at every site of $\mathcal P$
points to its partner. On this event, as long as $Y$ stays in
$\Lambda_m$, $X$ follows $Y$ with one two-step round trip inserted
on the first visit to each pair.
To verify this, suppose $X$ has used $k$ increments of $Y$, stands at
$Y_k=x\in\Lambda_m$, and the vertices whose arrows have been used form
a union of complete pairs. These statements hold at $k=0$. If the pair
containing $x$ is untouched, the two arrows send $X$ from $x$ to
$\phi(x)$ and back, consuming both arrows. If the pair was already
consumed, no such steps are needed. In either case $X$ next makes a
subsequent departure from $x$, uses the next increment of $Y$, and
reaches $Y_{k+1}$. This proves the description by induction while
$Y$ stays in the box. It applies also when $i\in B$, since the pairs
are disjoint even when both partners lie in $\Lambda_m$.

If $Y$ stays in $\Lambda_m$ through its $t$th step, the inserted
round trips stay in $\mathcal P$, and $X$ has taken at least $t$
physical steps by the time it uses the $t$th increment of $Y$.
Thus $N_t\le|\mathcal P|$. Independence of the arrows
and $Y$, together with \eqref{eq:box-survival}, now gives
\[
 \PP(N_{t}\le 2|\Lambda_m|)
       \ge (2d_1)^{-2|\Lambda_m|}e^{-Ct/m^2}.
\]
For large $r$ take the largest $m$ with $2(2m+1)^{d_2}\le r$.
Then $m$ is of order $r^{1/d_2}$, giving the lower bound.
\end{proof}

\begin{proof}[Proof of Proposition \ref{prop-rangeLaplacelog}]
Put $\kappa=d_2/(d_2+2)$ and $\Phi(t)=\EE e^{-\lambda N_{t}}$.
The two terms $r$ and $t/r^{2/d_2}$ in \eqref{eq:range-tail} have the
same order when $r$ is of order $t^\kappa$.
With $r=\lfloor t^\kappa\rfloor$, the lower bound in
\eqref{eq:range-tail} gives $\Phi(t)\ge e^{-\lambda r}\PP(N_{t}\le r)
\ge e^{-C_\lambda t^\kappa}$. For the upper bound take
$r=\lfloor\varepsilon t^\kappa\rfloor$. Splitting according to
whether $N_t\le r$ gives $\Phi(t)\le\PP(N_t\le r)+e^{-\lambda r}$.
Since $Cr-ct/r^{2/d_2}\le(C\varepsilon-c\varepsilon^{-2/d_2})t^\kappa$,
choosing $\varepsilon>0$ small bounds both terms by
$e^{-c_\lambda t^\kappa}$ for all sufficiently large $t$.
Taking negative logarithms and using $\EE e^{-\lambda R_t}=\Phi(t-1)$
proves the proposition.
\end{proof}

\section*{Acknowledgments}

Yuval Peres is supported by the National Natural Science Foundation of China under Grant Number W2531011. Shuo Qin is supported by the China Postdoctoral Science Foundation under Grant Numbers 2025M773086 and 2026T190815, and the National Natural Science Foundation of China under Grant Number 12271284. 

\emph{Use of AI tools.} An earlier version of Proposition \ref{lem:rectangle} was formulated through interaction with ChatGPT 6. In an earlier draft, the authors proved a weaker form of Proposition \ref{prop:window} by treating the candidate intervals in each column separately and estimating the bracket of the associated martingale. A review of that proof with ChatGPT led to the Fourier inversion argument used here. ChatGPT also contributed to the last-visit estimate, the explicit moment bounds in Section \ref{sec:column} and the block arguments in Lemma \ref{lem-range}. Claude Opus 5.5
(Anthropic) reviewed several versions of the manuscript and checked its
citations, drafted parts of the introduction, prepared
Figures \ref{fig:vertical-motion}, \ref{fig:column-identity}
and \ref{fig:pair-reversal}, revised Figures \ref{fig:model-rule}
and \ref{fig:model-simulation} and several captions, and contributed to the
exposition of the height martingale, the column occupation bound and the
defocusing argument. Claude Fable 5.2 audited one revision.
The authors checked all arguments and take full responsibility for the
paper.

\end{document}